\documentclass[a4paper]{article}
\usepackage{amsmath,amssymb,amsthm,color,graphicx,dsfont}
\usepackage[noPostScript]{diagrams}
\author{Michael Wibmer \thanks{Supported by the FWF (Project P16641)}}
\title{Gr\"{o}bner Bases for Families of Affine or Projective Schemes}

\newtheorem{lemma}{Lemma}
\newtheorem{theo}{Theorem}
\newtheorem{cor}{Corollary}
\newtheorem{defi}{Definition}
\theoremstyle{definition}
\newtheorem{ex}{Example}

\newcommand{\p}{\mathfrak{p}}
\newcommand{\q}{\mathfrak{q}}
\newcommand{\Spec}{\operatorname{Spec}}
\newcommand{\V}{\operatorname{V}}
\newcommand{\ida}{\mathfrak{a}}
\newcommand{\lt}{\operatorname{lt}}
\newcommand{\lc}{\operatorname{lc}}
\newcommand{\supp}{\operatorname{supp}}
\newcommand{\lm}{\operatorname{lm}}
\newcommand{\coef}{\operatorname{coef}}
\newcommand{\IY}{\mathcal{I}_Y}
\newcommand{\T}{\mathcal{T}}
\newcommand{\sigI}{\langle\sigma_\p(I)\rangle}
\newcommand{\G}{\mathcal{G}}

\begin{document}
\maketitle

\begin{abstract}
Let $I$ be an ideal of the polynomial ring $A[x]=A[x_1,\ldots,x_n]$ over the commutative, noetherian ring $A$.
Geometrically $I$ defines a family of affine schemes over $\Spec(A)$: For $\p\in\Spec(A)$, the fibre over $\p$ is the
closed subscheme of affine space over the residue field $k(\p)$, which is determined by the extension of $I$ under the
canonical map
$\sigma_\p:A[x]\rightarrow k(\p)[x]$. If $I$ is homogeneous there is an analogous projective setting, but again the
ideal defining the fibre is $\sigI$. For a chosen term order this ideal has a unique
reduced Gr\"{o}bner basis which is known to contain considerable geometric information about the fibre.
We study the behavior of this basis for varying $\p$ and prove the existence of a canonical decomposition of the base
space $\Spec(A)$ into finitely many locally closed subsets over which the reduced Gr\"{o}bner bases of the fibres can be
parametrized in a suitable way.
\end{abstract}
\section*{Introduction}
Let $A$ be a commutative, noetherian ring with identity and
$A[x]=A[x_1,\ldots ,x_n]$ the polynomial ring in the variables
$x_1,\ldots, x_n$ over $A$. We denote the residue field at $\p\in\Spec(A)$ by $k(\p)$. Geometrically an ideal
$I\subset A[x]$ defines a family of affine schemes over $\Spec(A)$: The canonical map $A\rightarrow A[x]/I$ gives
rise to a morphism of affine schemes
\[\varphi:\Spec(A[x]/I)\rightarrow\Spec(A).\]  For $\p\in\Spec(A)$ the fibre $\varphi^{-1}(\p)$ is the closed
subscheme of $\mathbb{A}^n_{k(\p)}=\Spec (k(\p)[x])$ determined by $\langle
\sigma_{\p}(I)\rangle$ where $\sigma_{\p}:A[x]\rightarrow
k(\p)[x]$ denotes the trivial extension of the canonical map
$A\rightarrow k(\p)$.

If $I$ is a homogeneous ideal we analogously obtain a family of projective schemes from
\[\varphi:\operatorname{Proj}(A[x]/I)\rightarrow\Spec(A).\]
The fibre $\varphi^{-1}(\p)$ is the closed subscheme of $\mathbb{P}^n_{k(\p)}=\operatorname{Proj}(k(\p)[x])$,
again determined by $\sigI$.

For a chosen term order we wish to study -- simultaneously
for all $\p\in\Spec(A)$ -- the unique reduced Gr\"{o}bner basis of $\sigI$. It is well known that such a Gr\"{o}bner
basis facilitates ``easy access'' to geometric information about the fibre $\varphi^{-1}(\p)$. It also seems reasonable
to
compare two fibres by ``comparing'' the corresponding Gr\"{o}bner bases. Of course we can compare
the leading terms, however it is not quite clear what comparing the Gr\"{o}bner bases should mean. We will make
this notion precise by introducing parametric sets. Rather vaguely a parametric set with respect to $I$ is a locally
closed subset $Y$ of $\Spec(A)$ such that over $Y$ the reduced Gr\"{o}bner bases of the fibres can be parameterized
in a suitable way. The main result of this article is to establish the existence and uniqueness of a canonical
decomposition of the base space $\Spec(A)$ into finitely many parametric sets.

Many concrete mathematical problems can be stated in the above described framework of families of affine or projective
schemes and to know the Gr\"{o}bner basis structure of the fibres may be the first step to their solution, if not
yet the solution itself. For example, if $A$ is a polynomial ring over some field, then we obtain the
case of algebraic systems with parameters, which is important for many ``real life" applications such as robotics or
electrical engineering (see e.g. \cite{cox-et}, chapter 6, and \cite{montes:loadflow}).
From a more theoretical point of view parametric sets are a tool to explore the geometry of families of affine or
projective schemes. Related theoretical applications range from efficient Gr\"{o}bner basis
computation (see e.g. \cite{arnold:modularalg} and \cite{pauer:lucky}) to cohomology (see \cite{walther:cohom}).

The naive hope that for a Gr\"{o}bner basis $G$ of $I$ the specialized Gr\"{o}bner basis $\sigma_\p(G)$ is a Gr\"{o}bner
basis of the specialized ideal $\sigI$ is in general not fulfilled. The behavior of Gr\"{o}bner bases under
specialization (or extension of scalars) has actually been studied by many authors, e.g.
\cite{bayer-et:extension}, \cite{ kalk:stability}, \cite{ aschenbr:reduction}, \cite{ gonzalez:homogeneous}
\cite{ gonzalez:hilbertstrat}, \cite{ assi},
\cite{ fortuna-et:degreduction}. In \cite{ aschenbr:reduction} the case of standard bases in the ring of formal
power series is treated. Relations to flatness are explored in \cite{ assi} and also in \cite{ bayer-et:extension}.
Articles focusing more on the fibres are \cite{ weispf:comprehensive}, \cite{ weispf:canonical},
\cite{ montes:anewalg} and \cite{ montes:canonical}. These last articles were written from a more computational point of
view, which led to a rather rash use of the word ``canonical''. So one main objective of the present article is
to establish a proper theoretical foundation for the underlying ideas of these articles.\\

\noindent The outline of the article is the following:
Section \ref{secparametric} \emph{(Parametric sets)} introduces the fundamental notion of parametric sets and their basic
properties. The main theorem of section \ref{seclucky} \emph{(Lucky primes and pseudo division)} is a
characterization of parametric sets in terms of lucky primes (see \cite{ grabe:lucky}). This theorem can also be
understood as giving the geometric meaning of luckiness.
Finally in section \ref{seccovers} \emph{(Gr\"{o}bner covers)} we achieve the main objective of the article by proving
existence and uniqueness of a canonical finite covering of $\Spec(A)$ with parametric subsets.

\subsection*{Preliminaries and notation}
A parametric subset $Y$ of $\Spec(A)$ facilitates an object which parameterizes the reduced Gr\"{o}bner bases
of $\sigI$ for $\p\in Y$.
To assure uniqueness of this object, which will be called the reduced Gr\"{o}bner basis of $I$ over $Y$ we have to
work with reduced schemes $(Y,\mathcal{O}_Y)$. In particular we would like to assume that our base
ring $A$ is reduced. This can be done without loss of generality:

Let $\operatorname{Nil}(A)$ denote the nilradical of $A$ and define $A'=A/\operatorname{Nil}(A)$. Then there is a natural
homeomorphism \begin{align*}\Spec(A) & \rightarrow \Spec(A')\\ \p & \mapsto \p '\end{align*}
and $k(\p)=k(\p')$. Moreover if $I'\subset A'[x]$ denotes the extension of $I$ under the canonical map
$A[x]\rightarrow A'[x]$ then $\sigI=\langle\sigma_{\p'}(I')\rangle$ for all $\p\in\Spec(A)$. \\

Throughout $A$ denotes a commutative, noetherian, reduced ring with
identity and $I$ an ideal of the polynomial ring $A[x]=A[x_1,\ldots,x_n]$.
For an $A$-module $M$ the localization at $\p\in\Spec(A)$ is denoted by $M_\p$ and
$k(\p)=A_\p/\p_\p$ is the residue field at $\p$. The map $\sigma_\p : A[x]\rightarrow k(\p)[x]$ denotes the
coefficientwise
extension of the canonical map $A\rightarrow k(\p)$.

We will only consider
reduced subschemes of $\Spec(A)$. So by a subscheme of $\Spec(A)$ we mean a locally closed subset $Y$ of
$\Spec(A)$ with the induced reduced subscheme structure $\mathcal{O}_Y$. Let $\ida$ be the radical ideal of
$A$ with $\overline{Y}=\V(\ida)$. (As usual $\V(\ida)\subset\Spec(A)$ denotes the closed set of all
prime ideals containing $\ida$.) For an open set $U$ of $Y$ we can explicitly describe $\mathcal{O}_Y(U)$ as the
set of all functions from $U$ into the disjoint union $\coprod (A/\ida)_\p$ which are locally fractions. We will
continuously identify $\Spec(A/\ida)$ with $\V(\ida)\subset\Spec(A)$.

The set of terms (i.e. powerproducts) is
denoted by $\T=\T(x_1,\ldots,x_n)$. Throughout we fix a term order
$<$ on $\T$. For a nonzero polynomial $P=\sum_{t\in \T}a_tt\in A[x]$ we
define
\begin{itemize}
\item the \emph{coefficient of $P$ at $t$} by $\coef(P,t)=a_t$,
\item the \emph{support of $P$} by $\supp(P)=\{t\in\T;\ a_t\neq 0\}$,
\item the \emph{leading term $\lt(P)$ of $P$} to be the maximal element of $\supp(P)$,
\item the \emph{leading coefficient of $P$} by $\lc(P)=\coef(P,\lt(P))$ and
\item the \emph{leading monomial of $P$} by $\lm(P)=\lc(P)\lt(P)$.
\end{itemize}
For $G\subset A[x]$ we set $\lt(G)=\{\lt(P);\ P\in G\smallsetminus\{0\}\}$ and similarly
$\lm(G)=\{\lm(P);\ P\in G\smallsetminus\{0\}\}$. A finite subset
$G$ of $I$ is called a Gr\"{o}bner basis of $I$ if $\langle\lm(G)\rangle=\langle\lm(I)\rangle$.
For $t\in\T$ we define the ideal of leading coefficients at $t$ by \[\lc(I,t)=\{\lc(P);\ P\in I \text{ with
} \lt(P)=t\}.\]
Note that $\lc(I,t)$ can conveniently be read off from a Gr\"{o}bner basis $G$ of $I$. In fact, $\lc(I,t)$ is generated
by \{$\lc(g);\ g\in G \text{ with } \lt(g)\text{ divides } t\}$. For a general reference for Gr\"{o}bner bases over
rings see \cite{ adams-et:intro}.\\

\noindent Before really getting started we look at some warm-up examples:

\begin{ex}
Let $k$ be a field and $A=k[u_1,u_2]$ the polynomial ring in the two parameters $u_1,u_2$. Consider the ideal
\[I=\big\langle (u_1^2-u_2)x, (u_2-1)y^2+u_1x\big\rangle\subset A[x,y].\] When faced with the task to
describe the Gr\"{o}bner basis structure of the fibres I guess most mathematicians would
come up with the following pictures:

\newpage
term order with $y^2>x$:
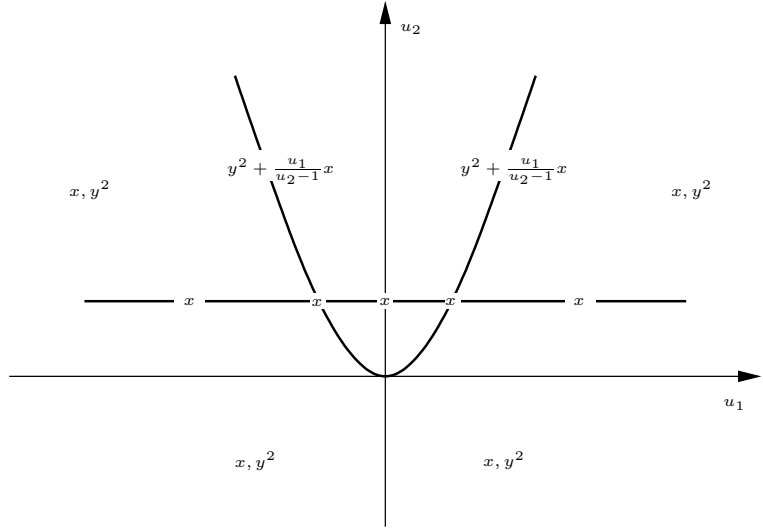
\begin{figure}[htbp]
\begin{center}

\input{ex1.pstex_t}

\end{center}
\end{figure}

term order with $x>y^2$:
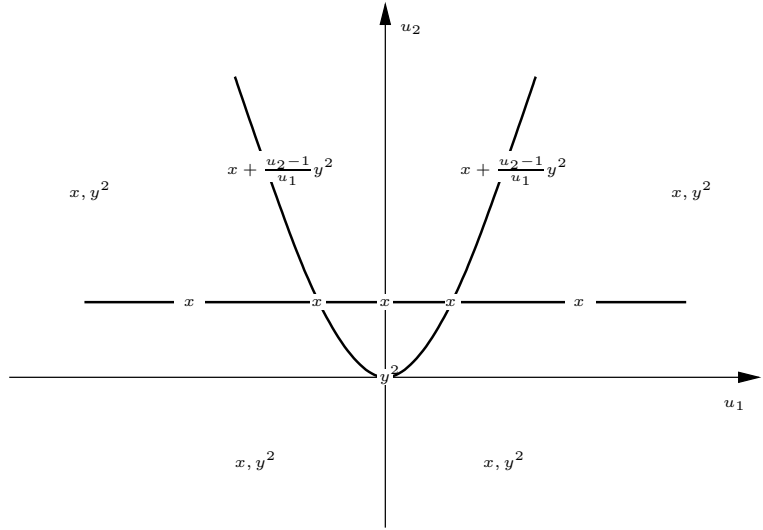
\begin{figure}[htbp]
\begin{center}

\input{ex1part2.pstex_t}

\end{center}
\end{figure}

The above pictures illustrate a decomposition of the base space $\mathbb{A}^2_k=\Spec(A)$ into
locally closed subsets. In short, the objective of this article is to find this decomposition in general.
\end{ex}

\begin{ex}\label{exlocal}
Let $k$ be an algebraically closed field and $A=k[u_1,u_2,u_3,u_4]$ the polynomial ring in the parameters
$u_1,u_2,u_3,u_4$. We consider the ideal \[I=\big\langle (u_2u_3-u_4u_1)x,\ u_1x^2+u_2x,\ u_3x^2+u_4x\big
\rangle\subset A[x].\]
(Here $x$ denotes just one variable.) Let $v=(v_1,v_2,v_3,v_4)\in k^4$ and
\[\p_v=\langle u_1-v_1,u_2-v_2,u_3-v_3,u_4-v_4\rangle.\] If $v_2v_3-v_4v_1$ is nonzero then the reduced
Gr\"{o}bner basis of $\langle\sigma_{\p_v}(I)\rangle$ is $x$. If $v_1$ and $v_3$ are zero and one of $v_2,v_4$ is
nonzero then the reduced Gr\"{o}bner basis of $\langle\sigma_{\p_v}(I)\rangle$ is also $x$. (In particular the set
of all $v\in k^4$ such that $\lt(\langle\sigma_{\p_v}(I)\rangle)$ is generated by $x$ is not locally closed.)
If $v$ lies in the quasi-affine variety
$Y=\V(\langle u_2u_3-u_4u_1\rangle)\smallsetminus\V(\langle u_1, u_3\rangle)$ then the reduced Gr\"{o}bner basis of
$\langle\sigma_{\p_v}(I)\rangle$
is given by $x^2+f(v)x$ where $f$ denotes the regular function on $Y$ defined by
\[f(v)=
\begin{cases} v_2/v_1 & \text{if } \ v_1\neq 0\\
v_4/v_3 & \text{if } \ v_3\neq 0 .
\end{cases}
\]
\end{ex}

\noindent The above example illustrates the ``local nature'' of the problem and suggests to work with sheaves and not
just with polynomials in $I$, as was common practice in \cite{ weispf:canonical} or \cite{ montes:canonical}.

\vspace{3mm}

Using the Buchberger algorithm it is relatively easy to see that the equivalence relation $\sim$ defined on
$\Spec(A)$ by $\p\sim\p'$ if $\lt(\sigI)=\lt(\langle\sigma_{\p'}(I)\rangle)$ has only finitely many equivalence classes
and that every
equivalence class is a constructible set. However there are reasons which militate against the
obvious approach to simply stratify the base space
$\Spec(A)$ with respect to the leading terms:

\begin{itemize}
\item The equivalence classes are indeed only constructible and not in general locally closed
(see example \ref{exlocal}).
\item Even if an equivalence class $Y$ is locally closed $\varphi$ may not be flat over $Y$.
\item Just because the function $\p\mapsto\lt(\sigI)$ is constant on $Y$ does not mean that the reduced Gr\"{o}bner
bases of the fibres depend on $\p\in Y$ in a ``continuous way''.
\end{itemize}

\noindent The following simple example illustrates the two latter points.

\begin{ex}\label{exhyp}
Let $k$ be a field and $A=k[u]$ the polynomial ring in one parameter $u$. Consider the ideal
$I=\langle u(ux-1), (ux-1)x\rangle\subset A[x]=k[u,x]$. The corresponding picture is on the next page.

\begin{figure}[h]
\begin{center}

\input{hyp.pstex_t}

\end{center}
\end{figure}

Geometrically the map $\varphi:\Spec(k[u,x]/I)\rightarrow\Spec(A)=\mathbb{A}_k^1$ is the projection onto
the $u$-axis.
For every point $\p\in\Spec(A)$ the leading terms of $\sigI$ are generated by $x$ but $\varphi$ is not flat: Let
$\mathfrak{P}\in\Spec(k[u,x]/I)$ be the point corresponding to the origin in $\mathbb{A}_k^2$ then
\[\mathcal{O}_{\Spec(A[x]/I),\mathfrak{P}}=
(A[x]/I)_{\mathfrak{P}}=k\] because $ux-1$ does not lie in $\mathfrak{P}$. For
$\p=\varphi(\mathfrak{P})=\langle u\rangle$ we have $\mathcal{O}_{\Spec(A),\p}=k[u]_{\langle u\rangle}$.
The map $k[u]_{\langle u\rangle}\rightarrow k$ induced by $\varphi$ is given by evaluation at the origin and is not flat.
Thus $\varphi$ is not flat at $\mathfrak{P}$.
\end{ex}

This example suggests that the above described problems may not appear in the projective setting. Indeed we
will see in section \ref{seccovers} that for homogeneous ideals the situation is as nice as could be hoped for, i.e.
the sets over which $\p\mapsto\lt(\sigI)$ is constant are parametric.

\section{Parametric sets}
\label{secparametric}
The idea of ``parameterizing Gr\"{o}bner bases'' can nicely be captured using sheaves.
For every subscheme $Y$ of $\Spec(A)$ we will define a quasi-coherent sheaf $\mathcal{I}_Y$ on $Y$, which
intuitively might be thought of as the restriction of $I$ to $Y$.

Let $Y$ be a locally closed subset of $\Spec(A)$ and $\ida\subset A$ the radical ideal such that $\overline{Y}=\V(\ida)$
and let $\overline{I}$ denote the extension of $I$ in $(A/\ida)[x]$. We define $\IY$ to be the restriction of the
quasi-coherent sheaf associated to the $A/\ida$ -module $\overline{I}$ on $\Spec(A/\ida)=\V(\ida)$ to $Y$. That is
\[\IY=\widetilde{\overline{I}}\Big|_Y.\]


More explicitly, for an open subset $U$ of $Y$ the $\mathcal{O}_Y(U)$-module $\mathcal{I}_Y(U)$ consists of all
functions $g$ from $U$ into the disjoint union $\coprod_{\p\in U}\overline{I}_{\p}$ which are locally fractions, i.e. for
every $\p\in U$ there exists an open neighborhood $U'$ of $\p$ in $U$ such that for
all $\q\in U'$ we have $g(\q)=\frac{P}{s}\in\overline{I}_\q$, where $P\in \overline{I}$ and
$s\in (A/\ida)\smallsetminus\q$ for all $\q\in U'$.


Since $A$ is noetherian, $\Spec(A)$ is a noetherian topological space and thus every open subset $U$ of $Y$ is
quasi-compact. This implies that we can consider $\mathcal{I}_Y(U)$ as an ideal of the polynomial ring
$\mathcal{O}_Y(U)[x]$.
(If $U$ was not quasi-compact we could not be sure that an element of $\mathcal{I}_Y(U)$ has finite support.)

Note that for $\p\in Y$ the stalk $\mathcal{I}_{Y,\p}=\overline{I}_\p$ is just the extension of $I$ under
$A[x]\rightarrow (A/\ida)_\p[x]$. Let $\mathfrak{m}_\p$ denote the unique maximal
ideal of $\mathcal{O}_{Y,\p}=(A/\ida)_\p$, then in analogy to the sequence
\[A\rightarrow\mathcal{O}_Y(Y)\rightarrow\mathcal{O}_{Y,\p}\rightarrow\mathcal{O}_{Y,\p}/\mathfrak{m}_\p=k(\p)\]
of natural maps we obtain natural maps
\[I\rightarrow\IY(Y)\rightarrow\mathcal{I}_{Y,\p}\rightarrow\sigI.\]
For $g\in\IY(Y)$ the image of $g$ in $\sigI$ is denoted by $\overline{g}^\p$.

Now we are prepared to give precise meaning to the intuitive idea of parameterizing Gr\"{o}bner bases: We are
looking for subschemes $Y$ of $\Spec(A)$ with the property that there exist global sections $g_1,\ldots,g_m\in\IY(Y)$
such that for all $\p\in Y$ their images $\overline{g_1}^{\p},\ldots,\overline{g_m}^{\p}$ form the unique reduced
Gr\"{o}bner basis of $\sigI$. We will need the following easy lemma.

\begin{lemma}\label{lembegin}
Let $Y$ be a subscheme of $\Spec(A)$ and $g,f\in\IY(Y)$. Then the set
\[\left\{\p\in Y;\ \overline{g}^\p=\overline{f}^\p\right\}\]
is a closed subset of $Y$ and $\overline{g}^\p=\overline{f}^\p$ for all $\p\in Y$ implies $g=f$.
\end{lemma}
\noindent Proof: It suffices to treat the case $f=0$. We can cover $Y$ with open sets $U_i$ such that
\[g(\p)=\frac{P}{s}\in\overline{I}_{\p}\] for $P\in\overline{I}\subset(A/\ida)[x]$ and
$s\in (A/\ida)\smallsetminus\p$ for all $\p\in U_i$. We have
\[\left\{\p\in Y;\ \overline{g}^\p=0\right\}\cap U_i=\left\{\p\in U_i;\ \coef(P,t)\in\p\text{ for all }
t\in\supp(P)\right\},\] which is a closed subset of $U_i$. Hence $\{\p\in Y;\ \overline{g}^\p=0\}$ is closed.

If we interpret $g$ as a polynomial with coefficients $c_t$ in $\mathcal{O}_Y(Y)$, then
$\overline{g}^{\p}=0$ is equivalent to saying that for all $t\in\T$ the image of $c_t$ in the stalk
$\mathcal{O}_{Y,\p}=(A/\ida)_{\p}$ lies in the maximal ideal $\mathfrak{m}_{\p}$ of $\mathcal{O}_{Y,\p}$.
Since this holds for all $\p\in Y$ and $Y$ is a reduced scheme we obtain $c_t=0\in\mathcal{O}_Y(Y)$. Hence $g=0$.\qed

\begin{theo}\label{theo1}
If $Y$ is a connected subscheme of $\Spec(A)$ and there exists a finite subset $G$ of $\IY(Y)$ such
that for all $\p\in Y$ the set $\overline{G}^\p=\{\overline{g}^{\p};\ g\in G\}$ is the reduced
Gr\"{o}bner basis of $\langle\sigma_{\p}(I)\rangle$, then $G$ is uniquely
determined and for each $g\in G$ the function $\p\mapsto\lt(\overline{g}^\p)$ is constant on $Y$. In particular, the
function $\p\mapsto\lt(\langle\sigma_\p(I)\rangle)$ is constant on $Y$.
\end{theo}
\noindent Proof: First we will show that for $g\in G$ and $t\in\T$ the set
\[W(t)=\{\p\in Y;\ \lt(\overline{g}^\p)=t\}\] is a closed subset of $Y$.
We can cover $Y$ with open sets $U_i$ such that
\[g(\p)=\frac{P}{s}\in\overline{I}_{\p}\text{ for all }
\p\in U_i.\] Here
$P\in\overline{I}\subset(A/\ida)[x]$ and $s\in(A/\ida)\smallsetminus\p$ for all $\p\in U_i$.

Let $\p\in Y$ and $\phi: (A/\ida)_\p\rightarrow(A/\ida)_\p/\mathfrak{m}_\p=k(\p)$ the canonical map. We will need that
$\phi(c/s)=1$ implies $c-s\in\p$ for $c\in A/\ida$ and $s\in(A/\ida)\smallsetminus\p$.
But $\phi(c/s)=1$ is equivalent to saying that there exists $c'\in\p$ and $s'\in(A/\ida)\smallsetminus\p$
such that \[\frac{c}{s}=1+\frac{c'}{s'}=\frac{s'+c'}{s'}.\] This implies the existence of an
$s''\in(A/\ida)\smallsetminus\p$ such that \[(cs'-s(s'+c'))s''=0\in\p.\]
Hence $cs'-ss'\in\p$ and therefore $c-s\in\p$.

Using the above result we see that for $\p\in U_i$
we have $\lt(\overline{g}^{\p})=t$ if and only if $\p$ contains $\coef(P,t')$ for $t'>t$ and
$\coef(P,t)-s$ (Use that $\overline{g}^{\p}$ is monic). Therefore $W(t)\cap U_i$ is a closed subset of $U_i$ and thus
$W(t)\subset Y$ is closed.

Since $\Spec(A)$ is a noetherian
topological space, a finite number of the $U_i$'s will do and therefore the function $\p\mapsto\lt(\overline{g}^\p)$
takes only finitely many values on $Y$. Consequently $Y$ is the disjoint
union of finitely many $W(t)$'s. By the connectedness assumption on $Y$ we can conclude that the function
$\p\mapsto\lt(\overline{g}^\p)$ is constant on $Y$.

Assume that, for $F\subset\IY(Y)$, $\overline{F}^\p$ is the reduced Gr\"{o}bner basis of
$\sigI$ for every $\p\in Y$. Then for $f\in F$ and a chosen $\p\in Y$ there exists a $g\in G$ such that
$\overline{f}^\p=\overline{g}^\p$. Since the leading terms of $\overline{f}^\p$ and $\overline{g}^\p$ are
independent of $\p$ this implies $\lt(\overline{f}^\p)=\lt(\overline{g}^\p)$ for all $\p\in Y$, but
as $\overline{F}^\p=\overline{G}^\p$ is the reduced Gr\"{o}bner basis we can conclude
$\overline{f}^\p=\overline{g}^\p$ for
all $\p\in Y$ and therefore $f=g\in G$ by lemma \ref{lembegin}.
\qed \\

The following example shows that both assertions of the above theorem may be false if $Y$ is not connected.
\begin{ex}
Let $Y=\{\p_1,\p_2\}$ where $\p_1$ and $\p_2$ are two distinct closed points of $\Spec(A)$. Note that
$\mathcal{O}_Y(Y)$ is just $k(\p_1)\times k(\p_2)$. For $j=1,2$ let $G_j$ denote the reduced Gr\"{o}bner bases
of $\langle\sigma_{\p_j}(I)\rangle$. Then for any subset $G$ of \[G_1\times G_2\subset
\langle\sigma_{\p_1}(I)\rangle\times\langle\sigma_{\p_2}(I)\rangle=\IY(Y)\] with the property that the
projections $G\rightarrow G_i$ are surjective
we have that $\overline{G}^\p$ is the reduced Gr\"{o}bner basis of $\sigI$ for every
$\p\in Y$.
\end{ex}

As we wish to have a definition suitable for all (not necessarily connected) subschemes of $\Spec(A)$ we simply demand
what we want.
\begin{defi} \label{parametric}
A locally closed subset $Y$ of $\Spec(A)$ is called \emph{parametric for Gr\"{o}bner bases with respect to $I$ (and $<$)}
if there
exists a finite subset $G$ of $\IY(Y)$ with the following properties:
\begin{enumerate}
 \item  $\overline{G}^{\p}$ is the reduced Gr\"{o}bner basis of $\langle\sigma_{\p}(I)\rangle$ for every $\p\in Y$.
\item For each $g\in G$ the function $\p\mapsto\lt(\overline{g}^{\p})$ is constant on $Y$.
\end{enumerate}
\end{defi}
Since the ideal $I\subset A[x]$ is clear from the context we usually omit the reference to $I$ and simply talk about
parametric subschemes of $\Spec(A)$.

\begin{theo}
Let $Y\subset\Spec(A)$ be parametric and $G$ a finite subset of $\IY(Y)$ satisfying the two conditions of the above
definition. Then $G$ is uniquely determined and the function $\p\mapsto\lt(\sigI)$ is constant on $Y$.
Furthermore every $g\in G$ is monic with $\lt(g)=\lt(\overline{g}^\p)$ for every $\p\in Y$.
\end{theo}
\noindent Proof: Because of condition (2) we can repeat the uniqueness proof as in the last paragraph of the proof of
theorem \ref{theo1}.

To show that every $g\in G$ is monic with $\lt(g)=\lt(\overline{g}^\p)$ observe that the coefficients of $g$ are
just elements of $\mathcal{O}_Y(Y)$. Since $(Y,\mathcal{O}_Y)$ is a reduced scheme every element of
$\mathcal{O}_Y(Y)$ is uniquely determined by its images in $k(\p)$ where $\p$ ranges over all of $Y$. \qed

\begin{defi}
Let $Y\subset\Spec(A)$ be parametric, then the uniquely determined subset $G=G_Y$ of $\IY(Y)$ of the above theorem is
called the \emph{reduced Gr\"{o}bner basis of $I$ over $Y$}. We define the \emph{leading terms of $Y$}, denoted by
$\lt(Y)$, to be the value of the constant function $\p\mapsto\lt(\sigI)$.
\end{defi}

To give the reader some idea where the journey is going we give the following definition at this early stage -- even
though we will not need it before section \ref{seccovers}.

\begin{defi}
A \emph{Gr\"{o}bner cover of $\Spec(A)$ with respect to $I$ (and $<$ )} is a finite set $\G$ of
pairs $(Y,G_Y)$ such that
$Y\subset\Spec(A)$ is parametric, $G_Y$ is the reduced Gr\"{o}bner basis of $I$ over $Y$ and
\[\bigcup_{(Y,G_Y)\in\G}Y=\Spec(A).\]
\end{defi}

\vspace{3mm}

\noindent Parametric sets are well behaved with respect to inclusion:

\begin{theo}\label{theorestrict}
Let $Y\subset\Spec(A)$ be parametric. Then every locally closed subset $Y'$ of $Y$ is parametric and the canonical map
$\IY(Y)\rightarrow\mathcal{I}_{Y'}(Y')$ maps the reduced Gr\"{o}bner basis of $I$ over $Y$ to the reduced Gr\"{o}bner
basis of $I$ over $Y'$.
\end{theo}

\noindent Proof: First of all let us construct the canonical map of the theorem: Assume $\overline{Y}=\V(\ida)$ and
$\overline{Y'}=\V(\ida')$ for radical ideals $\ida$ and $\ida'$ of $A$. Let $\overline{I}\subset (A/\ida)[x]$ and
$\overline{I}'\subset(A/\ida')[x]$ denote the corresponding extensions of $I$. As $\overline{Y'}\subset\overline{Y}$ we
have $\ida\subset \ida'$ and a canonical map $A/\ida\rightarrow A/\ida'$ which extends to $\varphi:\overline{I}
\rightarrow\overline{I}'$. Then for
$\p\in Y'\subset Y$ we have a canonical map \[\varphi_\p: \overline{I}_\p\rightarrow\overline{I}'_\p.\] Now an
element $g\in\IY(Y)$ gives rise to a function \[g':Y'\rightarrow\prod_{\p\in Y'}\overline{I}'_\p\] by $g'(\p)=\varphi_\p
(g(\p))$. One easily verifies that the map $\IY(Y)\rightarrow\mathcal{I}_{Y'}(Y'),\ g\mapsto g'$ is well
defined and a morphism.
For $\p\in Y'$ the commutative diagram

\begin{diagram}
\overline{I}_\p&       & \rTo  &       & \overline{I}'_\p\\
      & \rdTo &       & \ldTo & \\
      &       & \sigI &       &
\end{diagram}
gives rise to a commutative diagram
\begin{diagram}
\IY(Y)&       & \rTo  &       & \mathcal{I}_{Y'}(Y')\\
      & \rdTo &       & \ldTo & \\
      &       & \sigI &       &
\end{diagram}
From this the claim of the theorem follows. \qed

\vspace{3mm}

Next we will give a characterization of parametric sets in terms of monic ideals (see \cite{ pauer:lucky}).

%
%
%
%

\begin{defi}
An ideal $I\subset A[x]$ is called \emph{monic (with respect to $<$)} if
$\lc(I,t)\in\{\langle 0\rangle,\langle 1\rangle\}$ for all $t\in\T$. In other words: $I$ is monic if
for every $t\in\lt(I)$ there exists a
monic polynomial $P\in I$ with $\lt(P)=t$.
\end{defi}

There are quite a few definitions of reduced Gr\"{o}bner bases in the literature. We will use the one strictly
paralleling the field case.
\begin{defi}
A Gr\"{o}bner basis $G=\{g_1,\ldots,g_m\}$ of $I$ is called reduced if for $j=1\ldots,m$
\begin{itemize}
\item $g_j$ is monic and\item $\supp(g_j)\cap\lt(I)=\{\lt(g_j)\}$.
\end{itemize}
\end{defi}

With this definition not every ideal has a reduced Gr\"{o}bner basis, but as in the field case one
easily shows that
if it exists, it is unique and that $A[x]/I$ is a free $A$-module with basis $\T\smallsetminus\lt(I)$.
Concerning existence we have the following (cf. \cite{ pauer:lucky} and
\cite{ aschenbr:reduction}, theorem 2.11).
\begin{theo}
Let $I\subset A[x]$ be an ideal, then there exists a reduced Gr\"{o}bner basis of $I$ if and only if $I$ is monic.
\end{theo}
\noindent Proof: If there exists a reduced Gr\"{o}bner basis of $I$ then clearly $I$ is monic. Conversely if
$I$ is monic then we can choose monic polynomials $g_1,\ldots,g_m\in I$ such that $\lt(g_1),\ldots,\lt(g_m)$ is
the unique minimal generating set of $\lt(I)$. Now if we mutually reduce the $g_j$'s we end up with the desired
reduced Gr\"{o}bner basis of $I$.\qed

\vspace{3mm}

\noindent The connection to parametric subschemes is the following:
\begin{theo}\label{monic=parametric}
A subscheme $Y$ of $\Spec(A)$ is parametric if and only if $\IY(Y)\subset\mathcal{O}_Y(Y)[x]$ is monic, and in this
case the reduced
Gr\"{o}bner basis of $I$ over $Y$ is the reduced Gr\"{o}bner basis of $\IY(Y)$. In particular
$\lt(\IY(Y))=\lt(Y)$.
\end{theo}
\noindent Proof: Suppose that $Y$ is parametric and let $G\subset\IY(Y)$ denote the reduced Gr\"{o}bner basis of $I$
over $Y$. We will show that the leading term of every $f\in\IY(Y)$ is divisible by $\lt(g)$ for some $g\in G$. Since
$(Y,\mathcal{O}_Y)$ is a reduced scheme there exists a $\p\in Y$ such that the image of $\lc(f)\in\mathcal{O}_Y(Y)$
in $k(\p)$ is nonzero. For such a $\p$ we know that $\lt(f)=\lt(\overline{f}^\p)$ is divisible by
$\lt(\overline{g}^\p)=\lt(g)$ for some $g\in G$. Since the elements of $G$ are monic this shows that $\IY(Y)$ is monic.

Now suppose that $\IY(Y)$ is monic and let $G=\{g_1,\ldots,g_m\}$ denote the reduced Gr\"{o}bner basis of $\IY(Y)$.
For $f\in\IY(Y)$ the usual division (or reduction) algorithm shows that there exists a representation
\[f=f_1g_1+\cdots+f_mg_m\] such that for $i=1,\ldots,m$ we have $\lt(f_i)\lt(g_i)\leq\lt(f)$ and
\[\coef(f_i,t)\in\big\langle\coef(f,t');\ t'\geq t\lt(g_i)\big\rangle\text{ for all } t\in\T.\]
By the last condition we have
$\lt(\overline{f_i}^\p)\lt(\overline{g_i}^\p)\leq\lt(\overline{f}^\p)$ for every $\p\in Y$.
Because $\overline{f}^\p=\overline{f_1}^\p\overline{g_1}^\p+\cdots+\overline{f_m}^\p\overline{g_m}^\p$ this
shows that $\lt(\overline{f}^\p)$ is divisible by $\lt(\overline{g_i}^\p)$ for some $i\in\{1,\ldots,m\}$. Since
every element of $\sigI$ is of the form $\lambda\overline{f}^\p$ for $\lambda\in k(\p)$ and $f\in\IY(Y)$
we can conclude that $\overline{G}^\p$ is a Gr\"{o}bner basis of $\sigI$ for every $\p\in Y$. As $g\in G$ is monic
the function $\p\mapsto\lt(\overline{g}^\p)$ is clearly constant and since $G$ is reduced also $\overline{G}^\p$
is reduced. Thus we have shown that $Y$ is parametric and that $G$ is the reduced Gr\"{o}bner basis of
$I$ over $Y$. \qed

\vspace{3mm}

\noindent So the reduced Gr\"{o}bner basis $G$ of $I$ over $Y$ is actually a G\"{o}bner basis. In fact, by theorem
\ref{theorestrict}, $G|_{U}=\{g|_{U};\ g\in G\}$ is the reduced Gr\"{o}bner basis of
$\IY(U)\subset\mathcal{O}_Y(U)[x]$ for every open subset $U$ of $Y$.

\begin{cor}
$\Spec(A)$ is parametric with respect to $I$ if and only if $I$ is monic and in this case the reduced Gr\"{o}bner basis
of $I$ over $\Spec(A)$ is the reduced Gr\"{o}bner basis of $I$.
\end{cor}
\noindent Proof: This follows directly from the theorem because $\mathcal{I}_{\Spec(A)}(\Spec(A))=I$
(see \cite{ hartshorne}, chapter II, proposition 5.1). \qed

\vspace{3mm}

Next we will prove a local criterion for a locally closed subset of $\Spec(A)$ to be parametric. Using
this criterion we will then show that a family of affine or projective schemes over a parametric subset
of $\Spec(A)$ is flat. We need two easy lemmas.
\begin{lemma}\label{easylt}
Let $\p\in\Spec(A)$ and $f\in I_\p$. Then there exists $P\in I$ and $s\in A\smallsetminus\p$
such that
\[f=\frac{P}{s}\in I_\p\] and $\coef(P,t)=0$ whenever $\coef(f,t)=0$. In particular $\lt(P)=\lt(f)$.
\end{lemma}
\noindent Proof: By definition there exists $P\in I$ and $s\in A\smallsetminus\p$ such that $f=P/s\in I_\p$.
If $\coef(f,t)=\coef(P,t)/s\in A_\p$ is zero there exists an $s_t\in A\smallsetminus\p$ such that
$\coef(P,t)s_t=0$. If we multiply $P$ and $s$ by the product of all $s_t$'s where $t$ ranges over the support
of $P$ we obtain the desired representation of $f$. \qed

\begin{lemma}\label{easylt2}
Let $Y\subset\Spec(A)$ be locally closed and $\ida\subset A$ the radical ideal such that $\overline{Y}=\V(\ida)$.
Let $P\in\overline{I}\subset (A/\ida)[x]$. Then the leading term of the image of $P$ in $\IY(Y)$ equals the
leading term of $P$.
\end{lemma}
\noindent Proof: It suffices to show that there exists a $\p\in Y$ which does not contain $\lc(P)$. Assume the contrary,
then $Y$ is contained in the closed set \[W=\{\p\in\Spec(A/\ida);\ \lc(P)\in\p\}.\] But as $Y$ is dense in
$\V(\ida)=\Spec(A/\ida)$
we conclude that $W=\Spec(A/\ida)$ and thus $\lc(P)\in\p$ for all $\p\in\Spec(A/\ida)$. Because $\ida$ is
a radical ideal this yields the contradiction $\lc(P)=0$.

\begin{theo}\label{localcrit}
Let $Y\subset\Spec(A)$ be locally closed and $T'$ a set of terms such that $\T T'=T'$. Let $\ida\subset A$ denote
the radical ideal such that $\overline{Y}=V(\ida)$ and
$\overline{I}$ the extension of $I$ in $(A/\ida)[x]$. Then $Y$ is parametric with $\lt(Y)=T'$ if
and only if $\overline{I}_\p$ is monic
with $\lt(\overline{I}_\p)=T'$ for every $\p\in Y$.
\end{theo}
\noindent Proof: To show that
$\overline{I}_\p$ is monic with $\lt(\overline{I}_\p)=T'$ it suffices to prove $\lt(\overline{I}_\p)\subset T'$ because
this shows that the image of the reduced Gr\"{o}bner basis of $I$ over $Y$ in $\overline{I}_\p$ is the reduced
Gr\"{o}bner basis of
$\overline{I}_\p$. Let $P\in\overline{I}$ and $s\in (A/\ida)\smallsetminus\p$. By lemma \ref{easylt} we may assume
that the leading term of $P/s\in\overline{I}_\p$ equals the leading term of $P$. And by lemma \ref{easylt2}
the leading term of the image of $P$ in $\IY(Y)$ is the leading term of $P$. Thus $\lt(P/s)\in\lt(\IY(Y))=\lt(Y)=T'$.

For the converse direction let $T=\{t_1,\ldots,t_m\}$ denote the minimal generating set of $T'$. For $i=1,\ldots,m$
and $\p\in Y$ let $g_i(\p)$ denote the element of the reduced Gr\"{o}bner basis of $\overline{I}_\p$ with leading term
$t_i$. We want to show that $g_i$ defines an element of $\IY(Y)$. Let $\p\in Y$ and $P\in\overline{I}, s\in (A/\ida)
\smallsetminus\p$ such that $g_i(\p)=P/s\in\overline{I}_\p$. By lemma \ref{easylt} we may assume $\lt(P)=t_i$ and
$\coef(P,t)=0$ for $t\in T'\smallsetminus\{{t_i}\}$. Because
$g_i(\p)$ is monic there exists an $s'\in (A/\ida)\smallsetminus\p$ such that $(\lc(P)-s)s'=0$. The set
$U=\{\q\in Y; s,s'\notin\q\}$ is an open neighborhood of $\p$ in $Y$ and we have
$g_i(\q)=P/s\in\overline{I}_\q$ for all $\q\in U$ because $P/s\in\overline{I}_\q$ is monic with leading term $t_i$ and
$\supp(P/s)\cap T'=\{t_i\}$. Thus the $g_i$'s are elements of $\IY(Y)$.

For $f\in\IY(Y)$ there exists a $\p\in Y$ such that the image of $\lc(f)$ in $(A/\ida)_\p$ is nonzero. This implies
that the leading term of the image of $f$ in $\overline{I}_\p$ is the leading term of $f$ and thus we have
$\lt(f)\in\lt(\overline{I}_\p)=T'$.

Consequently $\lt(\IY(Y))=T'$ and because $g_i$ is monic with leading term $t_i$ for $i=1,\ldots,m$ by theorem
\ref{monic=parametric} we see that $Y$ is parametric. \qed

\vspace{3mm}

\noindent Recall that $\varphi$ denotes the map from $\Spec(A[x]/I)$ respectively
$\operatorname{Proj}(A[x]/I)$
to $\Spec(A)$.
\begin{cor}
If $Y\subset\Spec(A)$ is parametric then $\varphi$ is flat over $Y$, i.e. the map $\varphi^{-1}(Y)\rightarrow Y$
is a flat morphism.
\end{cor}
\noindent Proof: Let $\ida\subset A$ denote the radical ideal such that $\overline{Y}=\V(\ida)$ and let
$\overline{I}$ denote the extension of $I$ in $(A/\ida)[x]$.
The scheme structure on the set $\varphi^{-1}(Y)$ is given by identifying $\varphi^{-1}(Y)$ with
$X=\Spec(A[x]/I)\times_A Y$ or $X=\operatorname{Proj}(A[x]/I)\times_A Y$ respectively. Thus for
$\mathfrak{P}\in\varphi^{-1}(Y)$ the stalk $\mathcal{O}_{X,\mathfrak{P}}$
equals $((A/\ida)[x]/\overline{I})_\mathfrak{P}$ or $((A/\ida)[x]/\overline{I})_{(\mathfrak{P})}$.
(Here $((A/\ida)[x]/\overline{I})_{(\mathfrak{P})}$ denotes the elements of degree zero in the localized ring
$S^{-1}((A/\ida)[x]/\overline{I})$, where $S$ is the multiplicative system of all homogeneous elements of
$(A/\ida)[x]/\overline{I}$ which do not lie
in $\mathfrak{P}$.)
Let $\p=\varphi(\mathfrak{P})\in Y$. We have to show that $\varphi_\mathfrak{P}:\mathcal{O}_{Y,\p}\rightarrow
\mathcal{O}_{X,\mathfrak{P}}$ is flat.
In the affine case $\varphi_\mathfrak{P}$ can be factored:
\[\mathcal{O}_{Y,\p}=(A/\ida)_\p\rightarrow (A/\ida)_\p[x]/\overline{I}_\p=((A/\ida)[x]/\overline{I})_\p\rightarrow
((A/\ida)[x]/\overline{I})_\mathfrak{P}=\mathcal{O}_{X,\mathfrak{P}}.\]
By theorem \ref{localcrit} the ideal $\overline{I}_\p\subset (A/\ida)_\p[x]$ is monic and thus
$(A/\ida)_\p[x]/\overline{I}_\p$ is a free $(A/\ida)_\p$-module. In particular $(A/\ida)_\p[x]/\overline{I}_\p$
is a flat $(A/\ida)_\p$-module. 
Since ``localization is flat'' $((A/\ida)[x]/\overline{I})_\mathfrak{P}$ is a flat
$((A/\ida)[x]/\overline{I})_\p$-module. This completes the proof in the affine case.

In the projective case we know that $S^{-1}((A/\ida)[x]/\overline{I})$ is a flat $((A/\ida)[x]/\overline{I})_\p$-module
and therefore also a flat
$(A/\ida)_\p$-module. Since $((A/\ida)[x]/\overline{I})_{(\mathfrak{P})}$ is a direct summand of
$S^{-1}((A/\ida)[x]/\overline{I})$ also $((A/\ida)[x]/\overline{I})_{(\mathfrak{P})}$ is a flat $(A/\ida)_\p$-module.
\qed

\section{Lucky primes and pseudo division}
\label{seclucky}
Now it is time to introduce the concept of pseudo division (cf. \cite{ cox-et} and \cite{ montes:anewalg}).
This is basically just the usual division without fractions. The idea behind pseudo division already appeared
in the proof of theorem \ref{monic=parametric}.
\begin{defi}\label{pseudo division}
Let $f,g_1,\ldots,g_m\in A[x]$. A representation
\[cf=f_1g_1+\cdots+f_mg_m+r \] is called a \emph{pseudo division of $f$ modulo $g_1,\ldots,g_m$ (w.r.t. $<$)}
if the following assertions are satisfied:
\begin{itemize}
\item $f_1,\ldots,f_m,r\in A[x]$ and $c\in A$ is a product of leading coefficients of the $g_j$'s.
\item $\lt(f_j)\lt(g_j)\leq\lt(f)$ for $j=1,\ldots,m$.
\item No term in $\supp(r)$ is divisible by a leading term of the $g_j$'s.
\item $\coef(f_j,t)\in\big\langle\coef(f,t');\ t'\geq\lt(g_j)t\big\rangle$
for all $j\in\{1,\ldots,m\}$ and $t\in\T$.
\end{itemize}
\end{defi}

\noindent In this situation $r$ is called a remainder of $f$ after pseudo division modulo $g_1,\ldots,g_m$.
A pseudo division of $f$
modulo $g_1,\ldots,g_m$ can be obtained by successively applying pseudo reduction steps:

If there exists an element of the support of $f$ which is divisible by a leading term of any of
the $g_j$'s then choose $t\in\supp(f)$ maximal with this property. Then $t=t'\lt(g_j)$ holds for
some $j\in\{1,\ldots,m\}$ and $t'\in\T$. Now substitute $f$ by
\[\lc(g_j)f-\coef(f,t)t'g_j.\] By iterating this process and keeping track of the monomials used, we obtain
the desired representation.

The nice thing about pseudo reductions is that they are stable under specialization in the sense that
\[\lt(\overline{f_j})\lt(\overline{g_j})\leq\lt(\overline{f})\]
for $j=1,\ldots,m$. Here $\overline{g}$ denotes the coefficientwise reduction of $g\in A[x]$ modulo some ideal of $A$.
(This follows directly from the last assertion of the definition.)

Observe that $c$ may well be zero if $A$ is not an integral domain.
\begin{defi}
A prime ideal of $A$ is called lucky for $I$ if for every $t\in\lt(I)$ it does not contain $\lc(I,t)$.
\end{defi}
To my knowledge the expression ``lucky'' was coined by mathematicians working on modular algorithms to compute
Gr\"{o}bner
bases over $\mathbb{Q}$ (see \cite{ arnold:modularalg}, \cite{ pauer:lucky}, \cite{ grabe:lucky}).
Mod $p$-arithmetic avoids the phenomenon of coefficient growth but it is not a priori clear
which prime numbers $p$ can be used for lifting a Gr\"{o}bner basis over $\mathbb{Z}/\mathbb{Z}p$ to a Gr\"{o}bner basis
over $\mathbb{Q}$. So mathematicians must have considered themselves lucky when they picked a prime doing the job.

Let $T$ be the unique minimal generating set of $\lt(I)$. Because $\lc(I,t)\subset\lc(I,t')$
if $t$ divides $t'$, a prime $\p\in\Spec(A)$ is lucky for $I$ if and only if $\p$ does
not contain $\prod_{t\in T}\lc(I,t)$. In particular luckiness is an open condition.

\begin{defi}
The ideal \[J=J(I)=\sqrt{\prod_{t\in T}\lc(I,t)}\subset A\] is called the \emph{singular ideal of $I$ (with
respect to $<$)}.
\end{defi}
So a prime $\p\in\Spec(A)$ is unlucky (i.e. not lucky) for $I$ if and only if it is an element of the singular variety
$\V(J)$.

In \cite{ weispf:canonical} Weispfenning introduced another discriminant ideal which, however, can only be constructed if
$A$ is an integral domain. So for the time being assume that $A$ is an integral domain. In this case we can
consider the reduced Gr\"{o}bner basis $G$ of $I$ over the quotient field of $A$. For $g\in G$ the set
\[J_g=\{a\in A;\ ag\in I\}\] clearly is an ideal of $A$ and we can define Weispfenning's discriminant ideal by
\[J'=J'(I)=\sqrt{\prod_{g\in G}J_g}.\]

Clearly $J_g\subset\lc(I,\lt(g))$ always holds but the inclusion may be strict as illustrated by the following example.
\begin{ex}
Let $k$ be a field and $A=k[u_1,u_2]$ the polynomial ring in the parameters $u_1,u_2$. We consider the ideal
\[I=\langle u_1x+u_2,\ u_1y^2-1\rangle\subset A[x,y].\]
With respect to any term order the reduced Gr\"{o}bner basis of $I$ over the quotient field of $A$ is
\[G=\left\{x+\frac{u_2}{u_1},\ y^2-\frac{1}{u_1}\right\}.\]
But as $u_2y^2+x=y^2(u_1x+u_2)-x(u_1y^2-1)\in I$ we have with respect to any term order with $y^2>x$
\[J_{y^2-\frac{1}{u_1}}=\langle u_1\rangle\subsetneqq\langle u_1,u_2\rangle\subset\lc(I,y^2).\]
\end{ex}
However our discriminant ideal is not larger than Weispfenning's; in fact, they are the same.

\begin{theo}
In the above described situation we have $J=J'$.
\end{theo}
\noindent Proof: Let $I'$ denote the extension of $I$ in the polynomial ring over the quotient field of $A$.
First of all observe that $\lt(I)=\lt(I')$: As $I\subset I'$ the inclusion $\lt(I)\subset\lt(I')$ is clear. For the
other inclusion it suffices to notice that every $P\in I'$ is of the form $P=\frac{Q}{a}$ with $Q\in I$ and $a\in A$.

Let $G=\{g_1,\ldots,g_m\}$ denote the unique reduced Gr\"{o}bner basis of $I'$ over the quotient field
of $A$. Then as $\lt(I)=\lt(I')$ the unique minimal generating set $T$ of $\lt(I)$ equals
$\{\lt(g_1),\ldots,\lt(g_m)\}$. With the abbreviations $t_j=\lt(g_j)$ and $J_j=J_{g_j}$ for $j=1,\ldots,m$ we
may assume $t_1<\cdots <t_m$. We have to show \[\V\big(\lc(I,t_1)\cdots\lc(I,t_m)\big)=\V(J_1\cdots J_m).\]
As $J_j\subset\lc(I,t_j)$ for $j=1,\ldots,m$, the inclusion ``$\ \subset \ $'' is clear. For the other inclusion it
will suffice to show that for $j\in\{1,\ldots,m\}$ and $\p\in\Spec(A)$
\[J_j\subset\p \ \Rightarrow \ \lc(I,t_1)\cdots\lc(I,t_j)\subset\p.\]
We will prove this by contradiction. So assume $\lc(I,t_1)\cdots\lc(I,t_j)\nsubseteq\p$. Then we can find
$f_1,\ldots,f_j\in I$ with $\lt(f_i)=t_i$ and $\lc(f_i)\notin\p$ for $i=1,\ldots,j$. Pseudo reduction of
$f_j$ modulo $f_1,\ldots,f_{j-1}$ yields a polynomial $g\in I$ with $\lt(g)=t_j,\ \lc(g)\notin\p$ and no term
in $\supp(g)$ divisible by any $t_1,\ldots,t_{j-1}$. So no term in the support of $g-\lc(g)g_j\in I'$ is divisible
by any $t_1,\ldots,t_m$. Hence $\lc(g)g_j=g\in I$ and we conclude $\lc(g)\in J_j\subset\p$ (in contradiction to $\lc(g)
\notin\p$). \qed \\

The above theorem asserts that the concept of (in)essential specializations as introduced by
Weispfenning in \cite{ weispf:canonical} is equivalent to the older concept
of (un)lucky prime ideals. The advantage of the idea of luckiness is, of course, that it works for more general rings,
i.e. not only for integral domains. Observe that it is quite natural to work with rings which are not integral domains,
because even if you start with an integral domain (e.g. the polynomial ring over a field in some parameters), the
singular ideal $J$ will typically not be prime and therefore $A/J$ will not be an integral domain. The relevance of this
will become clear with the next theorem which gives a characterization of parametric subsets in terms of luckiness.

\begin{lemma}\label{ltequal}
Let $Y\subset\Spec(A)$ be parametric and $\ida\subset A$ the radical ideal such that $\overline{Y}=\V(\ida)$. If
$\overline{I}$ denotes the extension of $I$ in $(A/\ida)[x]$ then $\lt(Y)=\lt(\overline{I}).$
\end{lemma}
\noindent Proof: Let $t\in\lt(Y)$ and $\p\in Y$. From theorem \ref{localcrit} we know that
$\overline{I}_\p\subset (A/\ida)_\p[x]$ is monic with $\lt(\overline{I}_\p)=\lt(Y)$. Thus there exists
$P\in\overline{I}$ and $s\in (A/\ida)\smallsetminus\p$ such that the leading term of $P/s\in\overline{I}_\p$ equals
$t$. By lemma \ref{easylt} we may assume $t=\lt(P)\in\lt(\overline{I})$.

The inclusion $\lt(\overline{I})\subset\lt(Y)$ follows from lemma \ref{easylt2} and theorem \ref{monic=parametric}.
\qed
%
%
%

\vspace{3mm}

\noindent Now we are prepared to prove the main theorem of this section. It exhibits the
``geometric meaning'' of luckiness.

\begin{theo}\label{main1}
Let $Y$ be a locally closed subset of $\Spec(A)$ and $\ida\subset A$ the radical ideal such
that $\overline{Y}=\V(\ida)$. Denote by $\overline{I}$ the image of $I$ in $(A/\ida)[x]$. Then $Y$ is parametric for
Gr\"{o}bner bases with respect to $I$ if and only if \[Y\cap \V(J(\overline{I}))=\emptyset.\]
In other words: $Y$ is parametric if and only if every $\p\in Y$ is lucky for $\overline{I}$.
\end{theo}
\noindent Proof: Assume $Y$ is parametric and $\{g_1,\ldots,g_m\}\subset\IY(Y)$ is the reduced Gr\"{o}bner basis of
$I$ over $Y$. Then by lemma \ref{ltequal} the minimal generating set $T$ of
$\lt(\overline{I})$ equals $\{\lt(g_1),\ldots,\lt(g_m)\}$. Let $\p\in Y$ and $i\in\{1,\ldots,m\}$. By
lemma \ref{easylt} there exists $P_i\in\overline{I}$ with $\lt(P_i)=\lt(g_i(\p))$
and $s_i\in (A/\ida)\smallsetminus\p$ such that $g_i(\p)=P_i/s_i\in\overline{I}_\p$. Because
$\lt(P_i)=\lt(g_i(\p))=\lt(g_i)=\lt(\overline{g_i}^\p)$ we have $\lc(P_i)\notin\p$,
i.e. $\lc(\overline{I},\lt(P_i))\nsubseteq\p$.
Hence \[J(\overline{I})=\prod_{t\in T}\lc(\overline{I},t)\nsubseteq\p.\]
For the converse direction fix a $\p\in Y$ first and let $T=\{t_1,\ldots,t_m\}$ denote the minimal generating set
of $\lt(\overline{I})$. By assumption
\[\prod_{i=1}^m\lc(\overline{I},t_i)\nsubseteq\p.\]
Hence there exist polynomials $P_1,\ldots,P_m\in\overline{I}$ with $\lt(P_i)=t_i$ and
$\lc(P_i)\notin\p$. For $i=1,\ldots,m$ let $Q_i\in\overline{I}$ denote a remainder of $P_i$ after
pseudo division
modulo $\{P_1,\ldots,P_m\}\smallsetminus\{P_i\}$. Note that $\lt(Q_i)=\lt(P_i)=t_i$ and $\lc(Q_i)$ is a product of
leading coefficients of the $P_j$'s. Define \[U=\{\q\in Y;\ \lc(P_1)\cdots\lc(P_m)\notin\q\}.\] Then
$U$ is an open neighborhood of $\p\in Y$ and $Q_i/\lc(Q_i)$ defines an element of $\IY(U)$, which by
abuse of notation we again denote by $Q_i/\lc(Q_i)$.

We can repeat the above construction for any $\p'\in Y$ to obtain $U'$ and $Q_i'$ (analogously defined).
To obtain global sections $g_i\in\IY(Y)$ we have to show that
\[\frac{Q_i}{\lc(Q_i)}\bigg|_{U\cap U'}=\frac{Q_i'}{\lc(Q_i')}\bigg|_{U\cap U'}.\]
The leading term of \[\lc(Q_i')Q_i-\lc(Q_i)Q_i'\in\overline{I}\] is strictly smaller than $t_i$ and by construction no
term in the support of $\lc(Q_i')Q_i-\lc(Q_i)Q_i'$ is divisible by an element of
$\{t_1,\ldots,t_m\}\smallsetminus\{t_i\}$.
Thus $\lc(Q_i')Q_i-\lc(Q_i)Q_i'=0$ and we can glue together the sections $Q_i/\lc(Q_i)\in\IY(U)$ to obtain
global sections $g_i\in\IY(Y)$.

To show that $Y$ is parametric we will prove that $G=\{g_1,\ldots,g_m\}$ satisfies the conditions of definition
\ref{parametric}. Clearly $\lt(\overline{g_i}^\p)=t_i$ for every $\p\in Y$. So it remains to show that
$\overline{G}^\p$ is the reduced Gr\"{o}bner basis of $\sigI$ for every $\p\in Y$.
Let $\p\in Y$ and $P\in\overline{I}$. For a pseudo division (see definition \ref{pseudo division})
\[cP=P_1Q_1+\cdots+P_mQ_m+r\]
of $P$ modulo $Q_1,\ldots,Q_m$ we have $r\in\overline{I}$, but no term in the support of $r$ is divisible by an
element of $\{\lt(Q_1),\ldots,\lt(Q_m)\}=T$. Thus $r=0$ and
\[cP=P_1Q_1+\cdots+P_mQ_m.\] Let $\phi :(A/\ida)[x]\rightarrow k(\p)[x]$ denote the natural map then
\[\phi(c)\phi(P)=\phi(P_1)\phi(Q_1)+\cdots+\phi(P_m)\phi(Q_m)\]
and $\lt(\phi(P_i))\lt(\phi(Q_i))\leq\lt(\phi(P))$. Since $\lc(Q_i)\notin\p$ and
$c$ is a product of leading coefficients of the $Q_i$'s we know that $\phi(c),\phi(\lc(Q_1)),\ldots,\phi(\lc(Q_m))$ are
all nonzero. Consequently $\lt(\phi(P))$ is divisible by $\lt(\phi(Q_i))=t_i$ for some $i\in\{1,\ldots,m\}$.
Since every element of $\sigI$ is of the form $\lambda f$ for $\lambda\in k(\p)$ and
$f\in\phi(\overline{I})=\sigma_\p(I)$ this shows that $\lt(\sigI)$ is
generated by $T$ and so indeed $\overline{G}^\p$ is a Gr\"{o}bner basis of $\sigI$.
Clearly $\overline{g_i}^\p$ is monic and by
construction of the $Q_i$'s no term in the support of $\overline{g_i}^\p$ is divisible by an element of
$T\smallsetminus\{t_i\}$. Thus $\overline{G}^\p$ is the reduced Gr\"{o}bner basis of $\sigI$
and we are done. \qed

\begin{defi}
Let $Z$ be a closed subset of $\Spec(A)$ and $\ida\subset A$ the radical ideal such that $Z=\V(\ida)$. Let furthermore
$\overline{I}$ denote the extension of $I$ in $(A/\ida)[x]$. We define
\[Z_{gen}=Z\smallsetminus\V(J(\overline{I})).\]
\end{defi}

\begin{theo}\label{Zgen}
Let $Z\subset\Spec(A)$ be closed, $\ida\subset A$ the radical ideal such that $Z=\V(\ida)$ and $\overline{I}$
the extension of $I$ in $(A/\ida)[x]$. Then $Z_{gen}$ is parametric
with $\lt(Z_{gen})=\lt(\overline{I})$. Furthermore if $Y$ is an open subset of $Z$ such that
$Y$ is parametric with $\lt(Y)=\lt(\overline{I})$ then $Y\subset Z_{gen}$.

In other words: $Z_{gen}$ is the largest
open parametric subset of $Z$ with the same leading terms as $\overline{I}$.
\end{theo}
\noindent Proof:
Let $Y$ be an open subset of $Z$. First we will show that $\mathcal{I}_Z|_Y$ is canonically isomorphic to $\IY$.
Let $\ida'\subset A$ denote the radical ideal such that
$\overline{Y}=\V(\ida')$ and $\overline{I}'$ the extension of $I$ in $(A/\ida')[x]$.
Then $\ida\subset\ida'$ and the canonical map $A/\ida\rightarrow A/\ida'$ extends to
$\overline{I}\rightarrow\overline{I}'$ and further to $\phi:\overline{I}_\p\rightarrow\overline{I}'_\p$ for
$\p\in Y$. It suffices to show that $\phi$ is an isomorphism.

Clearly $\phi$ is surjective.
Let $P\in I$ and $s\in A\smallsetminus\p$ such that $\phi(\overline{P}/s)\in\overline{I}'_\p$ is zero.
This means that there exists $s'\in A\smallsetminus\p$ such that $\coef(s'P,t)\in\ida'$ for every $t\in\T$.
Let $\ida=\p_1\cap\cdots\cap\p_m$ be the (unique minimal) primary decomposition of the radical ideal $\ida$.
We may assume $\p_1,\ldots,\p_r\in Y$ and $\p_{r+1},\ldots,\p_{m}\notin Y$.
Note that $\p_i\notin Y$ implies $\V(\p_i)\cap Y=\emptyset$
because $Y$ is an open subset of $Z$. So in particular $\p_i\nsubseteq\p$ for $i=r+1,\ldots,m$. This means that
there exists an $s''\in\p_{r+1}\cap\cdots\cap\p_m\smallsetminus\p$. For $1\leq i\leq r$ we have
$\V(\p_i)\subset\overline{Y}=\V(\ida')$ and thus $\ida'\subset\p_i$.
Combining these results we see that every coefficient of $s''s'P$ lies in $\p_1\cap\cdots\cap\p_m=\ida$ and
thus $\overline{P}/s$ is zero in $\overline{I}_\p$. Consequently $\phi$ is injective.

An argument similar to the one above shows that for $\p\in Y$ the map $(A/\ida)_\p\rightarrow(A/\ida')_\p$ is an
isomorphism. Thus $\phi$ also preserves leading terms.

Now to show that $Z_{gen}$ is parametric
with the same leading terms as $\overline{I}$, just repeat the second part of the proof of theorem \ref{main1} (with
$Z_{gen}$ instead of $Y$) and use that $\mathcal{I}_Z(Z_{gen})$ is canonically isomorphic to
$\mathcal{I}_{Z_{gen}}(Z_{gen})$.

Now we additionally assume that $Y$ is a parametric subset of $Z$ with $\lt(Y)=\lt(\overline{I})$. Suppose
$Y\nsubseteq Z_{gen}$. Then there exists a $\p\in Y\smallsetminus Z_{gen}$. Let $T$ denote the minimal generating
set of $\lt(\overline{I})$. Since $\p\notin Z_{gen}=Z\smallsetminus\V(J(\overline{I}))$ there exists
a $t\in T$ such that $\lc(\overline{I},t)\subset\p$.

Since $Y$ is parametric with $\lt(Y)=\lt(\overline{I})$ we know from theorem \ref{localcrit} that
$\overline{I}'_\p$ is monic with $\lt(\overline{I}'_\p)=\lt(\overline{I})$. Using the isomorphism
$\phi:\overline{I}_\p\rightarrow\overline{I}'_\p$ we see that $\overline{I}_\p$ is monic with $\lt(\overline{I}_\p)=
\lt(\overline{I})$. Thus there exists $P\in\overline{I}$ and $s\in A\smallsetminus\p$ such that
$P/s\in\overline{I}_\p$ is monic with leading term $t$. By lemma \ref{easylt} we may assume $\lt(P)=t$.
Since $P/s$ is monic there exists $s'\in A\smallsetminus\p$ such that $(\lc(P)-s)s'=0$. Thus $\lc(P)\notin\p$
in contradiction to $\lc(\overline{I},t)\subset\p$. \qed

\vspace{3mm}

\noindent If we take $Z=\Spec(A)$ in the above theorem, then we see that the set of all lucky primes of $A$
($=\Spec(A)\smallsetminus\V(J(I))$) is the largest
open parametric subset of $\Spec(A)$ with the same leading terms as $I$. This more or less comes down to saying that
$J$ is the optimal discriminant ideal.\\

\noindent {\bf Caution:}
It is not true that $\p\in\Spec(A)$ is lucky for $I$ if and only if $\lt(I)=\lt(\sigI)$. We have seen above
that the ``only if'' direction is correct but the ``if'' direction is not true in general (see example \ref{exhyp}).
However it is true for homogeneous ideals as we will see in section \ref{subsechom}.

%
%
%
%
%
%
%
%

\vspace{3mm}

\noindent The following simple example illustrates that $Z_{gen}$ may well be the empty set.

\begin{ex}
Assume that $A$ is not an integral domain, then there exist $a,b\in A\smallsetminus\{0\}$ such that $ab=0$. If we take
$I$ to be the ideal of $A[x_1,x_2]$ generated by $ax_1$ and $bx_2$ then (with respect to
 any term order) $J(I)=\langle 0\rangle$
and so $\Spec(A)_{gen}=\emptyset$.
\end{ex}

However, this cannot happen if $Z$ is irreducible, because then $Z=\V(\ida)$ for some prime ideal $\ida$ of $A$ and
since
$A/\ida$ is an integral domain $J(\overline{I})$ is not the zero ideal and
thus $Z_{gen}$ is nonempty. In particular $Z_{gen}$ is dense in $Z$ and contains the generic point of $Z$.

The following examples have been included to convince the reader that the singular ideal $J$ is quite a reasonable
object.

\begin{ex}
Let $I\subset A[x]$ be the ideal generated by a square linear system

\[\begin{array}{c @{} c @{} c} P_1=b_{11}x_1+b_{12}x_2 & +\cdots+ &  b_{1n}x_n-c_1 \\
\vdots & \vdots  & \vdots \\
P_n=b_{n1}x_1+b_{n2}x_2&+\cdots+&b_{nn}x_n-c_n
\end{array}\]
and let \[B=\big(b_{ij}\big)_{1\leq i,j\leq n}\in A^{n\times n}\] denote the matrix of the system. Suppose
$\det=\det(B)\in A$ is not a zero divisor. Then the singular ideal $J$ of $I$ is independent of the chosen term order
and $\V(J)$ equals $\V(\det)$. In other words
$J=\sqrt{\langle \det\rangle}$.
\end{ex}
\noindent Proof: Let $B'\in A^{n\times n}$ denote the adjoint matrix of $B$. A classical linear algebra theorem
(see e.g. \cite{ lang:algebra}, chapter 8, $\S$ 4, proposition 8) asserts that
\begin{equation}\label{adjoint} B'B=BB'=\det\cdot\mathds{1},\end{equation}
where $\mathds{1}$ denotes the $n\times n$
identity matrix.

First we show that $1\notin\lt(I)$. Suppose the contrary. Let $A'$
denote the total ring of fractions of $A$, i.e. the localization
at the multiplicative subset of all nonzero divisors. Then we may
regard $A$ as a subring of $A'$. With the abbreviations
\[c=\left(\begin{array}{c} c_1 \\ \vdots \\ c_n
\end{array}\right)\text{ and } \xi=\frac{1}{\det}\cdot B'c\]
identity (\ref{adjoint}) shows that $\xi$ is a solution of our linear system.
Now $1\in\lt(I)$ simply means that there exist an $a\in A\smallsetminus\{0\}$ and $Q_1,\ldots,Q_n\in A[x]$ such that
\[Q_1P_1+\cdots+Q_nP_n=a.\] Evaluation at $\xi$ yields the contradiction $a=0$.

Identity (\ref{adjoint}) also shows that $\det$ lies in $\lc(I,x_i)$ for $i=1,\ldots,m$. Therefore $\det\in J$ and
$\V(J)\subset\V(\det)$. Now for the converse inclusion assume $\p\in\V(\det)$, i.e. $\det\in\p$. From theorem \ref{Zgen}
we know that for every $\q\in\Spec(A)\smallsetminus\V(J)$ the leading terms of $\langle\sigma_\q(I)\rangle$ are generated
by $x_1,\ldots,x_n$. But $\det\in\p$ implies that $\lt(\sigI)$ is not generated by $x_1,\ldots,x_n$ and consequently
$\p\in\V(J)$. \qed

%

\begin{ex}
Let $k$ be a field and $I'\subset k[x]=k[x_1,\ldots,x_n]$ a (homogeneous) ideal. For $1\leq i,j\leq n$ let $u_{ij}$
be additional indeterminates and abbreviate
\[ux=(u_{11}x_1+\cdots +u_{n1}x_n, \ldots\ldots  ,u_{1n}x_1+\cdots +u_{nn}x_n).\] Let $A$ be the polynomial ring over
$k$ in the $u_{ij}$'s and define
\[I=\langle P(ux);\ P\in I'\rangle\subset A[x].\]
Then 
the ideal of $k[x]$ generated by $\lt(\Spec(A)_{gen})$ is the generic initial ideal of $I'$, usually
denoted by $\operatorname{Gin}(I')$ (see e.g. \cite{ eisenbud:view} or \cite{ green:generic}).
\end{ex}

\begin{ex}
Suppose that $<$ is a graded order and $A$ is an integral domain, i.e. $\Spec(A)$ is irreducible.
Then $\Spec(A)_{gen}$ is a nonempty, open (and thus dense) subset of $\Spec(A)$ such that the function
\[\p\mapsto \text{affine Hilbert function of } \sigI\] is constant on $\Spec(A)_{gen}$. This is clear because
the affine Hilbert function of $\sigI$ is determined by $\lt(\sigI)$ (see \cite{ cox-et}, chapter 9, $\S$ 3,
proposition 4). Of course there is also an analogous ``projective'' statement.
\end{ex}


\section{Gr\"{o}bner covers}
\label{seccovers}
Now that we have (at least to some extent) explored the nature of parametric sets, it is time to see the complete
picture.

\begin{defi}
Let $L$ be a locally closed subset of $\Spec(A)$. A finite set $\G$ consisting of pairs $(Y,G_Y)$ with
$Y\subset\Spec(A)$ parametric and $G_Y$ the reduced Gr\"{o}bner bases of $I$ over $Y$ is called a
\emph{Gr\"{o}bner cover of $L$ with respect to $I$ (and $<$)} if
\[L=\bigcup_{Y\in\G}Y.\]
A Gr\"{o}bner cover $\G$ is called \emph{irreducible} if every $Y\in\G$ is irreducible.

A Gr\"{o}bner cover $\G$ of $L$ is called \emph{locally maximal} if for every $Y\in\G$ the following holds:
If $Y'\subset\Spec(A)$ is parametric with $Y'\subset L$ and $Y\subset Y'\subset\overline{Y}$ then $Y=Y'$.

A Gr\"{o}bner cover $\G$ is called \emph{small} if for every $Y\in\G$ we have
\[\overline{Y\smallsetminus\bigcup_{Y'\in\G\smallsetminus{\{Y\}}}Y'}=\overline{Y}.\]
\end{defi}

\vspace{3mm}

\noindent As already done in the above definition we write $Y\in\G$ instead of unhandy \mbox{$(Y,G_Y)\in\G$} and refer
to $Y$
as an element of $\G$. To say that a Gr\"{o}bner cover is small basically means that its elements are not
unnecessarily large. Our main interest, of course, is in Gr\"{o}bner covers of $\Spec(A)$ but
(with a view towards applications) it seems reasonable to also treat the relative case.

\begin{defi}
Let $L$ be a locally closed subset of $\Spec(A)$ and $G$ a finite subset of $I$. Then $G$ is called a \emph{comprehensive
Gr\"{o}bner basis of $I$ with respect to $L$ (and $<$)} if $\sigma_\p(G)=\{\sigma_\p(g);\ g\in G\}$ is a Gr\"{o}bner
basis
of $\sigI$ for every $\p\in L$.
\end{defi}
Comprehensive Gr\"{o}bner bases were introduced by Weispfenning in \cite{ weispf:comprehensive} and advanced in
\cite{ weispf:canonical}.
There is a rather obvious connection between Gr\"{o}bner covers of $L$ and comprehensive Gr\"{o}bner bases
of $I$ with respect to $L$, which we will now describe.

Let $\G$ be a Gr\"{o}bner cover of $L$. Choose a $Y\in\G$ and let $\ida\subset A$ be the radical ideal such that
$\overline{Y}=\V(\ida)$, furthermore let $\overline{I}$ denote the image of $I$ in $(A/\ida)[x]$.
Since $\Spec(A)$ is a noetherian topological space, $Y$ is
quasi-compact and so for every $g\in G_Y$ we can find finitely many open subsets $U_i$ of $Y$ which cover $Y$ and
have the following property: There exists a $P\in I$ and $s\in A/\ida$ such that
\[g(\p)=\frac{\overline{P}}{s}\in\overline{I}_\p \text{ for every } \p\in U_i.\]
Here $\overline{P}$ denotes the image of $P$ in $\overline{I}\subset(A/\ida)[x]$. Now taking together all
such $P$'s (for all $U_i$'s, all $g\in G_Y$ and all $Y\in\G$) we end up with a finite subset of $I$ which
clearly is a comprehensive Gr\"{o}bner basis of $I$ with respect to $L$.

%
%
%
%

\vspace{3mm}

The main theorem of this section
asserts that for every locally closed subset $L$ of $\Spec(A)$
there exists a unique irreducible, small and locally maximal Gr\"{o}bner cover of $L$. For the proof we will need a few
basic facts about constructible sets (cf. \cite{ hartshorne}).

\begin{defi}
Let $X$ be a topological space. A constructible subset of $X$ is a subset which belongs to the smallest
family $\mathfrak{F}$ of subsets such that
\begin{enumerate}
\item every open subset is in $\mathfrak{F}$,
\item a finite intersection of elements in $\mathfrak{F}$ is in $\mathfrak{F}$, and
\item the complement of an element in $\mathfrak{F}$ is in $\mathfrak{F}$.
\end{enumerate}
\end{defi}

\noindent One easily shows that the constructible sets of a topological space are exactly the finite unions
of locally closed sets.

\begin{lemma}\label{harts1}
Let $C$ be a constructible subset of $\Spec(A)$ and \[\overline{C}=Z_1\cup\cdots\cup Z_m\] the unique minimal
decomposition of $\overline{C}$ into irreducible and closed sets. Then for $j=1,\ldots,m$ there exists a nonempty open
subset of $Z_j$ contained in $C$.
\end{lemma}

\noindent Proof: A constructible set $C$ can be written as a finite union \[C=L_1\cup\cdots\cup L_{m'}\] of nonempty,
locally closed and irreducible sets $L_i$.
\[Z_1\cup\cdots\cup Z_m=\overline{C}=\overline{L_1}\cup\cdots\cup\overline{L_{m'}}\]
Fix a $j\in\{1,\ldots,m\}$. As $Z_j$ is irreducible there exists an $i\in\{1,\ldots,m'\}$ such that
$Z_j\subset\overline{L_i}$. Similarly, as $\overline{L_i}$ is irreducible there exist a $j'\in\{1,\ldots,m\}$ such that
$\overline{L_i}\subset Z_{j'}$. Hence \[Z_j\subset\overline{L_i}\subset Z_{j'}.\] This yields $j=j'$ and
$Z_j=\overline{L_i}$. So $L_i$ is a nonempty open subset of $Z_j$ contained in $C$. \qed

\begin{lemma}\label{harts2}
Let $L$ be a locally closed and irreducible subset of $\Spec(A)$. For a constructible subset $C$ of $\Spec(A)$ which
is contained in $L$ we have $\overline{C}=\overline{L}$ if and only if $C$ contains the generic point of $L$.
\end{lemma}
\noindent Proof: If $C$ contains the generic point $\p$ of $L$ we have $\overline{L}=
\overline{\{\p\}}\subset\overline{C}$. Hence by assumption $\overline{L}=\overline{C}$.

Conversely if $\overline{C}=\overline{L}$ by Lemma \ref{harts1} we know that there exists a nonempty open subset $U$ of
$\overline{L}$ contained in $C$. As $U\cap L$ is a nonempty open subset of $L$ we have \[\p\in U\cap L\subset C.\]
\qed

\begin{theo}\label{small}
Let $L\subset\Spec(A)$ be a locally closed set and $\G$ an irreducible Gr\"{o}bner cover of $L$. The following are
equivalent:
\begin{enumerate}
\item $\G$ is small.
\item Every $Y\in\G$ is the only element of $\G$ containing the generic point of $Y$.
\item For $Y,Y'\in\G$ with $Y\neq Y'$ and $Y\subset\overline{Y'}$ we have $Y\cap Y'=\emptyset$.
\end{enumerate}
\end{theo}
\noindent Proof: The equivalence of (1) and (2) follows from lemma \ref{harts2}.

For two distinct, locally closed and irreducible subsets $Y$ and $Y'$ of $\Spec(A)$ the generic point of $Y$ is contained
in $Y'$ if and only if $Y\subset \overline{Y'}$ and $Y\cap Y'\neq\emptyset$. Therefore (3) is equivalent to (2).\qed

\vspace{3mm}

\noindent Now we are prepared to prove the main theorem.

\begin{theo}
\label{main2}
Let $L$ be a locally closed subset of $\Spec(A)$. Then there exists exactly one irreducible, small and locally maximal
Gr\"{o}bner cover of $L$.
\end{theo}
\noindent Proof: First we will construct a Gr\"{o}bner cover $\G$ of $L$ and prove that it has the desired properties.
Then we will prove uniqueness. We construct $\G$ recursively:
\begin{center}
\fboxsep6pt
\fbox{\parbox{4in}{Set $C_1=L$ and $i=1$. \\
$(\star)$ Let
 \[\overline{C_i}=Z_{i1}\cup\cdots\cup Z_{im_i}\] be the unique minimal decomposition of
$\overline{C_i}$ into irreducible and closed sets. For $j=1,\ldots,m_i$ define
\[Y_{ij}=Z_{ij,gen}\cap\big(\text{union of all open subsets of } Z_{ij} \text{ contained in } L\big)\] and
\[C_{i+1}=C_i\smallsetminus (Y_{i1}\cup\cdots\cup Y_{im_i}).\]
If $C_{i+1}\neq\emptyset$ replace $i$ by $i+1$ and go to $(\star)$.}}
\end{center}
\vspace{2mm}
This yields a sequence of constructible sets $C_i$ with \[L=C_1\supset C_2\supset\cdots .\] To prove termination we will
show that the sequence \[\overline{C_1}\supset\overline{C_2}\supset\cdots\] is strictly decreasing. For $i\geq 1$ and
$j=1,\ldots,m_i$ there exists a nonempty open subset of $Z_{ij}$ contained
in $C_i\subset L$ by lemma \ref{harts1}. Hence $Y_{ij}$ is a nonempty open subset of $Z_{ij}$ contained in $L$.
\begin{align*}
\overline{C_{i+1}}& =\overline{C_i\smallsetminus(Y_{i1}\cup\cdots\cup Y_{im_i})}\subset
\overline{Z_{i1}\cup\cdots\cup Z_{im_i}\smallsetminus Y_{i1}\cup\cdots\cup Y_{im_i}} \\
& \subset \overline{(Z_{i1}\smallsetminus Y_{i1})\cup\cdots\cup(Z_{im_i}\smallsetminus Y_{im_i})}=
(Z_{i1}\smallsetminus Y_{i1})\cup\cdots\cup(Z_{im_i}\smallsetminus Y_{im_i})\\ & \subsetneqq Z_{i1}\cup\cdots\cup
Z_{im_i}=\overline{C_i}
\end{align*}
This shows that there exists a (minimal) $r\in\mathbb{N}$ such that $C_{r+1}=\emptyset$. Hence
\begin{align*}\emptyset& =C_{r+1}=C_r\smallsetminus (Y_{r1}\cup\cdots\cup Y_{rm_r})\\
& =C_{r-1}\smallsetminus (Y_{r-1,1}\cup\cdots\cup Y_{r-1,m_{r-1}}\cup Y_{r1}\cup\cdots\cup Y_{rm_r})=\cdots\\
& =C_1\smallsetminus (Y_{11}\cup\cdots\cup Y_{1m}\cup\cdots\cup Y_{r1}\cup\cdots\cup Y_{rm_r}).
\end{align*}
So we obtain \[L=C_1=Y_{11}\cup\cdots\cup Y_{1m}\cup\cdots\cup Y_{r1}\cup\cdots\cup Y_{rm_r}.\] As the $Y_{ij}$'s are
parametric by construction this shows that \[\G=\big\{(Y_{ij},G_{Y_{ij}})\ ;\ 1\leq i\leq r,\ 1\leq j\leq m_i\big\}\]
is a Gr\"{o}bner cover of $L$. It is clearly irreducible. Next we will show that $\G$ is locally maximal. So let
$Y\subset L$ be parametric with \[Y_{ij}\subset Y\subset\overline{Y_{ij}}=Z_{ij}.\] Then $Y$ is an open parametric
subset of $Z_{ij}$ and so by theorem \ref{Zgen} we have $Y\subset Z_{ij,gen}$. From the definition of
$Y_{ij}$ we obtain $Y\subset Y_{ij}$ and thus $Y=Y_{ij}$.\\

\noindent Now we will show that $\G$ is small. Let $Y_{ij},Y_{i'j'}\in\G$ with $(i,j)\neq(i',j')$.

\vspace{3mm}

We want to show that for $i\leq i'$ we have $Y_{ij}\nsubseteq\overline{Y_{i'j'}}$. Assume the contrary. Then
\[\overline{Y_{i'j'}}=Z_{i'j'}\subset\overline{C_{i'}}\subset\overline{C_i}=Z_{i1}\cup\cdots\cup Z_{im_i}.\]
Consequently there exists an $l\in\{1,\ldots,m_i\}$ such that $Z_{i'j'}\subset Z_{il}$. This yields
\[Z_{ij}=\overline{Y_{ij}}\subset\overline{Y_{i'j'}}=Z_{i'j'}\subset Z_{il}.\]
Therefore $j=l$ and $Z_{ij}=Z_{i'j'}$. For $i=i'$ this directly gives the contradiction $j=j'$. For $i< i'$ we have
\[Z_{ij}=Z_{i'j'}\subset\overline{C_{i'}}\subset\overline{C_{i+1}}\subset (Z_{i1}\smallsetminus Y_{i1})\cup\cdots\cup
(Z_{im_i}\smallsetminus Y_{im_i}).\]
Consequently $Z_{ij}\subset Z_{ij}\smallsetminus Y_{ij}$ and we obtain the contradiction $Y_{ij}=\emptyset$.

\vspace{3mm}

To prove that $\G$ is small it suffices, by theorem \ref{small}, to show that for $i>i'$ and
$Y_{ij}\subset\overline{Y_{i'j'}}$ we have $Y_{ij}\cap Y_{i'j'}=\emptyset$. Note that
$Y_{ij}\subset\overline{Y_{i'j'}}$ implies that $Z_{ij}\smallsetminus Y_{i'j'}$ is a closed subset of $\Spec(A)$. By
construction we have
\begin{equation}\label{Ci}
C_i=C_{i'}\smallsetminus \big(Y_{i'1}\cup\cdots\cup Y_{i'm_{i'}}\cup\cdots\cup Y_{i-1,1}\cup\cdots\cup
Y_{i-1,m_{i-1}}\big).
\end{equation}
For subsets $B,C,D$ of an arbitrary topological space with $D\subset C$ there is the trivial identity
\[\overline{\overline{B\smallsetminus C}\smallsetminus D}=\overline{B\smallsetminus C}.\]
Together with (\ref{Ci}) this yields
\begin{align*}
\overline{C_i} & =\overline{\overline{C_i}\smallsetminus Y_{i'j'}}=\overline{Z_{i1}\cup\cdots\cup Z_{im_i}\smallsetminus
Y_{i'j'}}\subset\overline{Z_{i1}\cup\cdots\cup(Z_{ij}\smallsetminus Y_{i'j'})\cup\cdots\cup Z_{im_i}}\\ & =
Z_{i1}\cup\cdots\cup(Z_{ij}\smallsetminus Y_{i'j'})\cup\cdots\cup Z_{im_i}\subset Z_{i1}\cup\cdots\cup Z_{im_i}=
\overline{C_i}.
\end{align*}
Therefore \[Z_{i1}\cup\cdots\cup Z_{im_i}=Z_{i1}\cup\cdots\cup(Z_{ij}\smallsetminus Y_{i'j'})\cup\cdots\cup Z_{im_i}\]
and $Z_{ij}\subset Z_{ij}\smallsetminus Y_{i'j'}$. Thus $Y_{ij}\cap Y_{i'j'}=\emptyset$.

\vspace{3mm}

So far we have shown that $\G$ is an irreducible, small and locally maximal Gr\"{o}bner cover of $L$. It remains to prove
uniqueness. Assume $\G'$ is another irreducible, small and locally maximal Gr\"{o}bner cover of $L$. First we will show
$\G\subset\G'$. More precisely we will show, by induction on $i=1,\ldots,r$, that $Y_{i1},\ldots,Y_{im_i}\in\G'$.
We denote
the generic point of $Y_{ij}$ by $\p_{ij}$.

\vspace{3mm}

\noindent First assume $i=1$. Let $j\in\{1,\ldots,m_1\}$. As \[\bigcup_{Y\in\G}Y=L=\bigcup_{Y'\in\G'}Y'\] there exists a
$Y'_{1j}\in\G'$ such that $\p_{1j}\in Y'_{1j}$. We want to show $Y_{1j}=Y_{1j}'$. As $Y_{1j}'$ is irreducible
and $\overline{Y_{1j}'}\subset\overline{L}=Z_{11}\cup\cdots\cup Z_{1m_1}$ there exist a $j'\in\{1,\ldots,m_1\}$
such that $\overline{Y_{1j}'}\subset Z_{1j'}$. Together with $\p_{1j}\in Y'_{1j}$ this gives
\[Z_{1j}\subset\overline{Y_{1j}'}\subset Z_{1j'}.\] Therefore $j=j'$ and $\overline{Y_{1j}'}=Z_{1j}$. Thus $Y_{1j}'$ is
an open subset of $Z_{1j}$ contained in $L$ and by theorem \ref{Zgen}
$Y_{1j}'\subset Z_{1j,gen}$. So by definition of $Y_{1j}$ we have $Y_{1j}'\subset Y_{1j}$. Since $\G'$ is locally maximal
we obtain $Y_{1j}=Y_{1j}'\in\G'$.

Now we do the induction step. Suppose \[Y_{11},\ldots,Y_{1m_1},\ldots,Y_{i-1,1},\ldots,Y_{i-1,m_{i-1}}\in\G'.\]
We have to show $Y_{i1},\ldots,Y_{im_i}\in\G'$. For $j\in\{1,\ldots,m_i\}$ there exists a $Y_{ij}'\in\G'$ such that
$\p_{ij}\in Y_{ij}'$. Using that $\G'$ is small and the induction hypothesis we obtain
\[\overline{Y_{ij}'}=\overline{Y_{ij}'\smallsetminus\bigcup_{Y'\in\G'\smallsetminus\{Y_{ij}'\}}Y'}\subset
\overline{L\smallsetminus\bigcup_{{1\leq i'\leq i-1 \atop 1\leq j'\leq m_{i'}}}Y_{i'j'}}=\overline{C_i}=
Z_{i1}\cup\cdots\cup Z_{im_i}.\] Hence there exists a $j'\in\{1,\ldots,m_i\}$ such that $\overline{Y_{ij}'}\subset
Z_{ij'}$. Together with $\p_{ij}\in Y'_{ij}$ this gives
\[Z_{ij}\subset\overline{Y_{ij}'}\subset Z_{ij'}.\] Therefore $j=j'$ and $\overline{Y_{ij}'}=Z_{ij}$. Since
$\G'$ is locally maximal a similar argument as in the case $i=1$ above proves $Y_{ij}=Y_{ij}'\in\G'$.
Thus we have shown $\G\subset\G'$.

Assume this is a proper inclusion. Then there exist a $Y'\in \G'$ such that $Y'\notin\G$ and therefore
\[\overline{Y'}=\overline{Y'\smallsetminus\bigcup_{Y\in\G'\smallsetminus\{Y'\}}Y}\subset
\overline{Y'\smallsetminus\bigcup_{Y\in\G}Y}=\overline{Y'\smallsetminus L}=\emptyset.\] This is a contradiction
as, by definition, the empty set is not irreducible. \qed

\begin{defi}
Let $L$ be a locally closed subset of $\Spec(A)$. The uniquely determined irreducible, small and locally maximal
Gr\"{o}bner cover of $L$ is called the \emph{canonical irreducible Gr\"{o}bner cover of $L$ (with respect to
$I$ and $<$)}.
\end{defi}

In \cite{ weispf:canonical} Weispfenning gave a rather ad hoc kind of construction for what he called
canonical Gr\"{o}bner systems. This construction bears some analogy with the existence proof of the above
theorem, however there are some differences between the concept of canonical Gr\"{o}bner systems and the concept of
canonical irreducible Gr\"{o}bner covers. For example, the canonical Gr\"{o}bner system may contain redundant
elements. The persistent reader is invited to verify this with the example $A=k[u_1,u_2]$ and $I=\langle u_1u_2,
u_1x^2+x\rangle$. (The point is simply that if $\Spec(A)=Z_1\cup\cdots\cup Z_m$ is the decomposition of $\Spec(A)$
into irreducible closed sets, then it may happen that the singular part of $Z_i$ ($= Z_i\smallsetminus Z_{i, gen}$)
is contained in some $Z_{j,gen}$.)

Note that theorem \ref{main2} implies that the equivalence relation on $\Spec(A)$, given by comparing the leading terms
of $\sigI$, has only finitely many equivalence classes and that every equivalence class is a constructible set.
Indeed example \ref{exlocal} and \ref{exnmin} show that these equivalence classes are only constructible and not locally
closed.
The following example illustrates that the canonical irreducible Gr\"{o}bner cover may be not of minimal
cardinality among the irreducible Gr\"{o}bner covers.

\begin{ex}\label{exnmin}
Let $k$ be a field and $A=k[u_1,u_2]$ the polynomial ring in the two parameters $u_1,u_2$. We consider the ideal
\[I=\langle u_1x, (u_2^2-1)x^2+x\rangle\subset A[x].\]
(Here $x$ denotes just one variable.) Obviously $J=J(I)=\langle u_1\rangle$ and the affine plane without the
$u_2$-axis has generic Gr\"{o}bner basis $x$, i.e. $Y_1=\mathbb{A}^2_{gen}=\Spec(A)\smallsetminus\V(u_1)$ and
$x\in\mathcal{I}_{Y_1}(Y_1)=I_{u_1}$ ($=$ localization of $I$ at $\{1,u_1,u_1^2,\ldots\}$) is the reduced Gr\"{o}bner
basis of $I$ over $Y_1$. By factoring mod $J=\langle u_1\rangle$ and identifying $A/J$ with $k[u_2]$ we obtain
\[\overline{I}=\langle(u_2^2-1)x^2+x\rangle\subset k[u_2][x].\]
On the $u_2$-axis the generic Gr\"{o}bner basis is $x^2+\frac{1}{u_2^2-1}x$, i.e.
\[J(\overline{I})=\langle u_2^2-1\rangle=\langle u_2+1\rangle\cap\langle u_2-1\rangle,\]
$Y_2=\V(u_1)_{gen}
=\V(u_1)\smallsetminus\V(u_2^2-1)$ and $x^2+\frac{1}{u_2^2-1}x\in\mathcal{I}_{Y_2}(Y_2)=\overline{I}_{u_2^2-1}$ is
the reduced Gr\"{o}bner bases of $I$ over $Y_2$. Finally over the two closed points $Y_3=\langle u_1, u_2-1\rangle$ and
$Y_4=\langle u_1,u_2+1\rangle$ we have the reduced Gr\"{o}bner basis $x$ again. To summarize
\[\G=\big\{(Y_1,\{x\}), (Y_2,\{x^2+{\textstyle \frac{1}{u_2^2-1}}x\}), (Y_3,\{x\}), (Y_4,\{x\})\big\}\] is the canonical
irreducible Gr\"{o}bner cover of $\mathbb{A}^2=\Spec(A)$.

Let $f\in k[u_1,u_2]$ be an irreducible polynomial such that $f(0,u_2)=u_2^2-1$ (e.g. $f=u_1+u_2^2-1$).
Then there exist $h\in A=k[u_1,u_2]$ such that $f=hu_1+u_2^2-1$, thus $fx^2+x=(hx)(u_1x)+(u_2^2-1)x^2+x\in I$.
Therefore the extension of $I$ in $(A/\langle f\rangle)[x]$ is just $\langle x\rangle$ and $\V(f)$ is parametric
with reduced Gr\"{o}bner basis $x$. Consequently
\[\G'=\{(Y_1,\{x\}),(Y_2,\{x^2+{\textstyle \frac{1}{u_2^2-1}}x\}), (\V(f),x)\}\]
is an irreducible Gr\"{o}bner cover of $\mathbb{A}^2$ with smaller cardinality than the canonical irreducible
Gr\"{o}bner cover. However, choosing an irreducible Gr\"{o}bner cover of $\Spec(A)$ with minimal cardinality
in a canonical way is as impossible as choosing a curve which meets the $u_2$-axes only in $(0,-1)$
and $(0,1)$ in a canonical way.

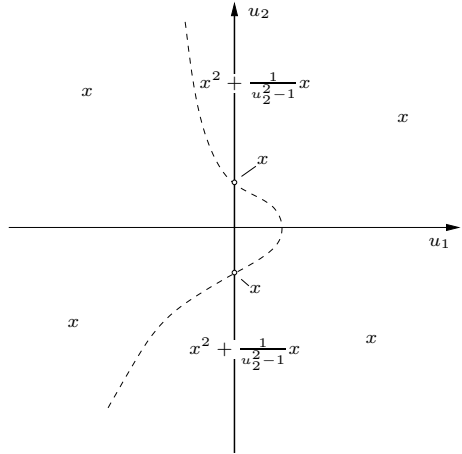
\begin{figure}[htbp]
\begin{center}

\input{nmin.pstex_t}

\end{center}
\end{figure}
\end{ex}

\noindent The above example can also be used to show that a parametric subset of $\Spec(A)$ need not be contained in
a maximal parametric subset.

\subsection{The projective case}\label{subsechom}

In the projective setting, i.e. if $I$ is a homogeneous ideal the situation is considerably nicer than in the
affine setting. It actually is as nice as it
can be hoped for: The equivalence classes of the equivalence relation $\sim$ defined on $\Spec(A)$ by $\p\sim\p'$ if
$\lt(\sigI)=\lt(\langle\sigma_{\p'}(I)\rangle)$ are parametric. (In particular they are locally closed.) The key to
the proof is the following lemma which is not true for arbitrary ideals (cf. example \ref{exhyp} and \ref{exnmin}).
The equivalence of (1) and (2) has already
been proved for $A=\mathbb{Z}$ in \cite{ arnold:modularalg} (theorem 5.13).

\begin{lemma}\label{lemhomok}
Let $I\subset A[x]$ be a homogeneous ideal and $\p\in\Spec(A)$. Then the following
assertions are equivalent:
\begin{enumerate}
\item $\p$ is lucky for $I$.
\item $\lt(\sigI)=\lt(I)$.
\item $\lt(\sigI)\supset\lt(I)$.
\end{enumerate}
\end{lemma}
\noindent Proof: We have already seen that (1) implies (2) in theorem \ref{Zgen}.
So we only have to show that (3) implies (1):

Assume that $\p\in\Spec(A)$ is unlucky for $I$. Then there exists $t\in\lt(I)$ such that
$\lc(I,t)\subset\p$. We may assume that $t$ is maximal in its degree, i.e. for every $t'\in\lt(I)$ with
$\deg(t')=\deg(t)$ and $\lc(I,t')\subset\p$ we have $t'\leq t$.
Since $t\in\lt(I)\subset\lt(\sigI)$ there
exists $P\in I$ such that $\lt(\sigma_\p(P))=t$. Because $I$ is homogeneous we may assume that $P$ is homogeneous and
thus $\deg(P)=\deg(t)$. We can also assume that $\lt(P)$ is minimal, i.e. for $P'\in I$ with $\lt(\sigma_\p(P'))=t$
we have $\lt(P')\geq\lt(P)$.

Because $\lc(I,t)\subset\p$ we have $\lt(P)>t$. By the maximality of $t$ we conclude $\lc(I,\lt(P))\nsubseteq\p$.
Thus there exists $Q\in I$ with $\lt(Q)=\lt(P)$ and $\lc(Q)\notin\p$. Set \[P'=\lc(Q)P-\lc(P)Q.\]
Then for $t'>t$ we have $\coef(P',t')\in\p$ because $\coef(P,t'),\lc(P)\in\p$.
On the other hand $\coef(P',t)$ does not lie in $\p$ because $\lc(Q),\coef(P,t)\notin\p$. Therefore
$\lt(\sigma_\p(P'))=t$ but as $\lt(P')<\lt(P)$ this contradicts the minimality of $P'$.\qed

\vspace{3mm}

Note that if $I\subset A[x]$ is an arbitrary ideal and $\p\in\Spec(A)$ is unlucky for $I$ then we can say virtually
nothing about the relation between $\lt(\sigI)$ and $\lt(I)$. We may have $\lt(\sigI)\subsetneqq \lt(I)$.
(This for example happens if $I$ is a monomial ideal.) Or we may have $\lt(\sigI)\supsetneqq\lt(I)$. (This for example
happens if $I$ is generated by a single polynomial $P=\sum_{i=1}^m a_it_i$ such that $t_i$ divides $t_{i+1}$ and
the $a_i$'s generate the unit ideal in $A$.) It may also
happen that $\lt(\sigI)$ and $\lt(I)$ are incomparable, i.e. there does not hold any inclusion relation between them.
Finally it may actually happen that $\lt(\sigI)$ equals $\lt(I)$ (see example \ref{exhyp}).

By the above lemma, we at least know that $\lt(I)$ is not contained in $\lt(\sigI)$ if $I$ is homogeneous and $\p$
unlucky for $I$.

\begin{theo}
Let $I\subset A[x]$ be a homogeneous ideal and $L\subset\Spec(A)$ locally closed. Then the equivalence classes of
the equivalence relation $\sim$ defined on $L$ by $\p\sim\p'$ if $\lt(\sigI)=\lt(\langle\sigma_{\p'}(I)\rangle)$
are parametric with respect to $I$.
\end{theo}
\noindent Proof: By theorem \ref{theorestrict} every locally closed subset of a parametric subset is parametric. Thus
we may assume $L=\Spec(A)$. Let $Y\subset\Spec(A)$ be an equivalence class and $T'\subset\T$ such that
$\lt(\sigI)=T'$ for all $\p\in Y$. From theorem \ref{main2} we already know that $Y$ is a constructible subset
of $\Spec(A)$. Let $Z$ be the closure of $Y$ and $\ida\subset A$ the radical ideal such that $\overline{Y}=Z=\V(\ida)$.
As usual $\overline{I}$ denotes the extension of $I$ in $(A/\ida)[x]$. To apply lemma \ref{lemhomok}
we have to show $\lt(\overline{I})=T'$. Let
\[Z=Z_1\cup\cdots\cup Z_m\] be the unique minimal decomposition of $Z$ into irreducible and closed subsets.
For $i=1,\ldots,m$ let $\ida_i\subset A$ denote the radical ideal such that $Z_i=\V(\ida_i)$ and $\overline{I}_i$ the
extension of $I$ in $(A/\ida)[x]$. By
lemma \ref{harts1} the intersection $Z_{i,gen}\cap Y$ is nonempty. Therefore by theorem \ref{Zgen} we have
$\lt(\overline{I}_i)=\lt(Z_{i,gen})=T'$.

Now let $P\in\overline{I}$. If for every $i\in\{1,\ldots,m\}$ the leading term of the image of $P$ in $\overline{I}_i$
is strictly smaller than the leading term of $P$, then the leading coefficient of $P$ must
lie in the intersection of all the $\ida_i$'s which is zero mod $\ida$. Thus there exists $i\in\{1,\ldots,m\}$
such that $\lt(P)\in\lt(\overline{I}_i)=T'$. Consequently $\lt(\overline{I})\subset T'$.

For the converse direction let $t\in T'=\lt(\overline{I}_1)$. There exists $P\in\overline{I}$ such that the
leading term of the image of $P$ in $\overline{I}_1$ is $t$. This means $\coef(P,t')\in\ida_1$ for $t'>t$ and
$\coef(P,t)\notin\ida_1$. The $\ida_i$'s constitute the minimal primary decomposition of $\ida$ and so
we can find $c\in\ida_2\cap\cdots\cap\ida_m\smallsetminus\ida_1$.
For $t'>t$ the coefficient of $cP$ at $t'$ lies in the intersection of all the $\ida_i$'s and thus equals zero.
On the other hand $\coef(cP,t)$ does not lie in $\ida_1$ and therefore $\lt(cP)=t$. Consequently
$t\in\lt(\overline{I})$.

By definition $Y$ is the set of all primes $\p\in Z$ such that $\lt(\sigI)$ equals $T'=\lt(\overline{I})$. Thus, by lemma
\ref{lemhomok}, $Y$ is the set of all lucky primes of $\overline{I}$, i.e. $Y=Z_{gen}$ which is parametric by proposition
\ref{main1}.\qed

\vspace{3mm}

It is now obvious how to define the canonical Gr\"{o}bner cover in the projective case:

\begin{defi}
Let $I$ be a homogeneous ideal of $A[x]$ and $L$ a locally closed subset of $\Spec(A)$. The Gr\"{o}bner cover
corresponding to the stratification of $L$ with respect to the leading terms of $\sigI$ is called the
\emph{canonical Gr\"{o}bner cover of $L$ with respect to $I$ (and $<$).}
\end{defi}

\subsection*{Conclusion and open questions}
We have introduced two concepts for studying the geometry of fibres:
parametric sets and Gr\"{o}bner covers. It seems possible to generalize these notions to
more general (i.e. not necessarily affine) base schemes.

Clearly one of the main reasons for the success of Gr\"{o}bner bases in the last decades has been the fact that
in many cases they could actually be computed. The focus of this article was not on algorithms, but of course
an efficient implementation of an algorithm to compute Gr\"{o}bner covers is desirable. The existence proof
for the canonical irreducible Gr\"{o}bner cover is in principle constructive, but an
algorithm for the computation of the canonical irreducible Gr\"{o}bner cover would necessarily involve successive
primary decompositions and thus would be of modest practical value.
The obvious solution is to skip irreducibility. For the projective case we have the canonical Gr\"{o}bner cover at hand
and it suggests itself to exploit this for the affine case by a process of homogenizing and dehomogenizing.

The problem of determining the Gr\"{o}bner basis structure
of the fibres has already been considered from an algorithmic point of view (see \cite{ montes:canonical},
\cite{ montes:anewalg}, \cite{ weispf:canonical}, \cite{ weispf:comprehensive}). Most notably Antonio Montes released
an implementation in Maple (see \texttt{ http://www-ma2.upc.edu/$\sim$montes}) for the important case where $A$ is the
polynomial ring over $\mathbb{Q}$. In fact, the output of his algorithm BUILDTREE can be interpreted as a Gr\"{o}bner
cover, but a drawback is that you cannot say a priori which Gr\"{o}bner cover the algorithm will compute, furthermore
the result depends on a term order on the parameters.
\subsection*{Acknowledgement}
I am grateful for the financial support from the Austrian Science Fund FWF (Project P16641) and the
caring backup of the project leader Kurt Girstmair. Furthermore I would like to give credit to
RISC and RICAM for supporting my stay at the ``Special Semester
on G\"{o}bner Bases and Related Methods'' in Linz and to all the people who discussed the topic of this article
with me.

%

\bibliographystyle{plain}
\bibliography{bibdata}

\bigskip

\noindent Michael Wibmer\\
Institute of Mathematics\\
University Innsbruck\\
Technikerstr. 25\\
A-6020 Innsbruck, Austria\\

\noindent e-mail: \texttt{csac8728@uibk.ac.at}\\
homepage: \texttt{http://relationproject.at.tt}
\end{document}